\documentclass[a4paper,10pt]{article}
\usepackage[utf8]{inputenc}
\usepackage{amsmath,amssymb}
\title{Some remarks on maximal rank}
\author{Tove Dahn (Lund University)}

\begin{document}
\maketitle

\section{Introduction}

\subsection{Continuity}
Assume $\Psi : U \rightarrow U^{\bot}$
projective with $\dim U^{c} = \dim U^{\bot}$, $U + U^{c}=I$. 
When $(I + \Psi) U \sim I$, we have $U \sim \Sigma c_{j} \Psi^{j}$. Thus, $d \Psi^{j} \in \mathcal{E}^{(0)'}$
implies $d U \in \mathcal{E}^{(0)'}$. Assume $\Phi$ a class of movements. Let $\Phi_{0}=\{ U \in \Phi \quad d U \in \mathcal{E}^{(0)'} \}$.
Define $\Phi_{ac}=\{ U \in \Phi \quad d U = 0 \quad U^{\bot}=0 \}$ and
$\Phi_{N}=\{ U U^{\bot} = U^{\bot} U \}$ (implies $U I = I U$ given U absolute continuous (a.c.))
Note the difference between parametrices and fundamental solutions, that is the parametrices ; 
$\mathcal{D}' \rightarrow \mathcal{D}^{F'}$ (\cite{Dahn15}).
Given $\Phi$ transitive through $\Psi : U \rightarrow U^{\bot}$
with $U^{2} g \simeq U g$, $g \in (I)$, where (I) is an ideal, such that $U + U^{\blacktriangle}=I$ over $(I)$, where $U^{\blacktriangle}=
(I-U) =-\Psi U \in \Phi$, that is $U^{\blacktriangle}$ projective on $(I)^{c}$. More precisely, Assume $R$ 
the restriction to $(J) \subset (I)$,
then $U$ is projective on $(J)$ iff $U R = R U$ ($R^{\bot} U=U R^{\bot}=0$).

\emph{Connected components:} Assume that the components are domains for constant character of movement . 
We consider the closure in a domain for constant dimension, that is we assume the movement is relatively 
closed, for instance
$d U=\alpha d V$, with $dV$ analytic and $\alpha$ constant close to the boundary. We assume $\mathcal{G}$ is defined
by movements with $dU$ of bounded variation (BV). Assume in particular, that $U_{1}$ is translation and 
$U_{2}$ rotation and that change of character is through
$d U_{2}=\alpha d U_{1}$, with $\alpha$ constant close to the boundary, (cf. \cite{Nishino68}). Note 
that if $d V \sim P d x$, for a polynomial P, we have that if $d V=0$ on $\Omega$, then $\Omega$ is locally algebraic. 
Note that constant dimension for $(U,U^{\bot})$ does not imply that the movement does not change character.

A \emph{very regular boundary} implies existence of $U$ analytic at the boundary. Assume further existence $\forall \Omega_{j}$ 
of $U_{j} \in \mathcal{G}$, such that $\Omega_{j}$ is generated by $U_{j}$, with $d U_{j}=\rho d U_{1}$ and $\rho$ regular or constant.
Note that if $U^{2}$ and U are of the same character, this does not imply
that U is projective. More precisely, given that $U^{\bot}(f)=\int f d U^{\bot}=0$ on $\Omega$ with restriction 
$R=R_{\Omega}$ and given $U^{\bot}R \equiv 0$ implies $(I-U)R f=0$, that is $\mid I-U \mid \leq \mid U^{\bot} \mid$
we have $(U + U^{\bot}) R=R$. 
Connected sets can be given by $U \gamma \rightarrow \gamma$,
when $U \rightarrow I$ continuously. Sufficient, given $F(\gamma)$ analytic and
$U F(\gamma)=F({}^{t} U \gamma)$, with ${}^{t} U \gamma(\zeta)=\gamma(\zeta_{T})$,
is that $\gamma$ is analytically dependent of $T$.

The movement is determined by $\eta d x - \xi d y=0$, that is $d U=0$ with the condition $Y/X \simeq \eta / \xi$.
$\frac{d}{d t} U=0$ is a differential operator in a real parameter, that is hypoelliptic.
Consider for this reason $\frac{d}{d t_{1}} + i \frac{d}{d t_{2}}(U,U^{\bot})=0$,
that is $\frac{d U}{d t_{1}} - \frac{d U^{\bot}}{d t_{2}}=0$.
It is sufficient that $U$ is independent of $t_{2}$ and
$U^{\bot}$ is independent of $t_{1}$, for instance $t_{2}=t_{1}^{*}$.
Assume $\Omega_{j} = \cup \Omega_{jk}$ and $d W_{k}=0$ on $\Omega_{jk}$, we then have $\Sigma d W_{k}=0$
on $\Omega_{j}$, as long as the sum is finite. Given $d W_{k}=0$ implies $W_{k}=I$ on $\Omega_{jk}$, we have that
$W=d W=0$ are isolated. 

A domain is \emph{simply connected}, if every simple Jordan curve, that divides the domain in two 
parts, can be continuously deformed to a point. Assume $\gamma \subset nbhd \infty$, such that 
$\gamma \rightarrow \infty$, as $t \rightarrow \infty$ continuously, then $\gamma$ defines on $t \geq 0$, a 
connected domain. Given $1/\gamma \rightarrow 0$, as $t \rightarrow \pm \infty$ for instance, 
$\gamma$ defines a simply connected  domain. Note $1/(\gamma + \frac{1}{\gamma}$) defines a simply 
connected domain, where $\gamma \rightarrow 0 , \pm \infty$. Consider in particular $\gamma$ a reduced 
polynomial. We consider $1 / (\gamma + 1/\gamma)$ as closed. 

\subsection{Factorization}

Given a short exact sequence
$0 \rightarrow T \rightarrow H \rightarrow S \rightarrow 0$, we have $H=T \bigotimes S (\simeq T \bigotimes \widehat{T})$.
In particular we have for $U \rightarrow I \rightarrow U^{\bot}$, given
$U^{\bot} \rightarrow 0$ regularly, that $I = U \bigoplus U^{\bot}$
that is U projective.
Assume $V$ surjective on (I) and $VW \simeq \mathcal{F}I$, where $\mathcal{F}$ is the Fourier transform. 
In particular $\mathcal{F} U I \simeq U^{\bot} V W$ and $\mathcal{F} I U \simeq V  W U$, where we assume
$\widehat{U T}=U^{\bot} \widehat{T}$.
The condition $(V W)^{\bot}= V W$ (spiral), gives that $(V W U)^{\bot}=U^{\bot} V W$.
Note that given ${}^{t} V=U_{1}$ and ${}^{t} W=U_{2}$, we have that ${}^{t}(V W)=V W$. 

$F(x,y,z,p,q)=\lambda$ can be solved through $\Phi,\Psi$ mutually independent and in involution with each other 
and with F. A complete solution to $F=\lambda$ is given by $\Phi=\mu,\Psi=\nu$, with parameters $\mu,\nu$. 
Conversely we have that when $\Phi,\Psi$ are solutions to $F=0$, completed through $\mu,\nu$, 
they generate solutions to $F=\lambda$. 
For instance when $\Phi$ is a normal system and $\Psi$ a not normal system.
When $\Gamma$ is defined by $d \Phi=0$, the closure can be for instance 
$w(p,q)=w(a,b)=0$ algebraic, completed to
an analytic zero, alternatively we consider pseudo convex domains, for instance such that $d \Phi^{\bot}=0$. 

Consider (Legendre) $L^{*}(x)=<x,x^{*}>$ and
$L^{\diamondsuit}$ the scalar product corresponding to harmonic conjugation (\cite{AhlforsSario60}). 
We then have that over an involutive set
$\{ L^{*},L^{\diamondsuit} \}=0$, that is $\xi^{*} \xi + \eta^{*} \eta=0$, so that to determine the movement, 
it is sufficient to
consider $L^{*}(x)=const$ (\cite{Lie91}). As long as the continuation $L^{*} \rightarrow L^{\bot}$ preserves the 
Legendre relation, the movement can be related to translation, $d U = \rho d U_{1}$ close to the boundary. 
Note (\cite{Lie91}) that given $L^{*},L^{\diamondsuit}$ are mutually independent and mutually in involution
and with F, given solutions to $L^{*}=b,L^{\diamondsuit}=c$, we have a complete solution to $F=a$,
where $a,b,c$ are constants. The solution is dependent of the parameters $b,c$.

\section{Comparable movements}
\subsection{Majorisation principle}
Given f almost periodic (a.p) (n=1), we have $M(f(x))=\lim_{T \rightarrow \infty} M_{T}(f(x))=\lim_{T \rightarrow \infty} \frac{1}{T} \int_{0}^{T} f(x) d x$. Thus,
if $F'=f$, where f a.p., we have that F a.p. iff $M_{T}(f) \rightarrow 0$, $T \rightarrow \infty$ or when F a.p. in $C^{\infty}$,
$M(\frac{d F}{d x}) \in \dot{B}$.
Assume $f=g/h$, where $h,g$ a.p.
 Given $h \neq 0$, we have that g/h a.p., that is $M(f/g) \rightarrow 0$ and we have a \emph{majorisation principle}. 
 The boundary to a pseudo convex domain is cylindrical (\cite{Oka60}), We will here limit ourselves to pseudo  convex domains.
 \newtheorem{def1}{Definition}[section]
 \begin{def1}
  Assume $\mid U^{\bot} g \mid \sim \mid w_{U} \mid \mid g \mid$ and in the same manner for $V,w_{V}$.
  We say that $U$ is weaker than $V$, $U \prec V$, if $\mid w_{U} \mid \leq \mid w_{V} \mid$, in $\infty$.
  We say that $U$ is strictly weaker than $V$, if $U \prec \prec V$, if $w_{U} / w_{V} \rightarrow 0$ in $\infty$. 
 \end{def1}

 We assume $I \prec U$, that is $U$ preserves compact sub level sets. Consider $V_{\lambda}=\{ \widehat{f} \leq \lambda \}$
 and assume $d U^{\bot}=0$ on $V_{\lambda}$. Given $\widehat{f}$ a.p. on $V_{\lambda}$ there is a sequence $\gamma_{j}$
 with a limit such that $U^{\bot} \widehat{f}=U f$. Given $d \mu=(g/h) d t$, that is $h d \mu=g d t$, 
 we can for instance choose h so that \emph{$\Omega$ is normal relative $h d \mu$}. Assume $u_{1},u_{2}$
 coordinates relative translation and rotation.
 Given $\Omega \rightarrow \tilde{\Omega}=\{ (u_{1},u_{2}) \}$ 
 contains a spiral $u_{1}=u_{2}$,
 by choosing $h$ so that $h=0$ on $u_{1},u_{2}$, we have that $h d \mu$ defines a normal surface $\tilde{\Omega}$,
 given that $h d \mu$ is analytic outside $u_{1}=u_{2}$.

Assume $d U=\beta d U_{1}$,
where $\beta$ regular and bounded, then a maximum principle for translates can be continued to $U$. 
When $U$ monotonous (increasing), we have that $U^{\bot}$ monotonous (decreasing). When $U_{S}=U(s,t)$ is a spiral,
with $U=U^{\bot}$, we have that 
that $U$ is monotonous in $(s,t)$.

\subsection{Mean convergence}
Assume $\tau f = \rho f$ with $\rho$ regular. When $\frac{d \rho}{d t}=0$ implies
$\tau=1$ is a point, then $\rho$ is invertible. Assume $d U = \rho d V$ and $d V=\vartheta d U$. 
Given $\mid \log \rho \mid \in L^{1}$
we have that $\rho, \frac{1}{\rho} \in L^{1}$, that is we can choose
$\rho \simeq 1/\vartheta$. 
Relatively compact translation implies uniformly continuous limit. It is sufficient that
$\mid d U \mid \rightarrow \mid d I \mid$, which is implied by $\mid d U^{\blacktriangle} \mid \rightarrow 0$.
Note that $d U$ relatively compact implies existence of $\gamma_{j} \rightarrow \gamma_{0}$, that is for instance
a maximum-principle.

\emph{Assume $U$ projective in the mean, that is $M(U f + U^{\bot} f)=M(f)$}. 
For instance $U f=\alpha * f$, 
where we assume $<U f,g>=<f, \beta * g>$ with $\alpha,\beta \in \mathcal{D}$. Further, assume $d U + d V= d U_{1} \sim d x$.  Given $V f$ 
harmonic, we have that $M(V f)=0$ implies $V f=0$.
Where $d U$ a reduced BV measure, $M_{U}(f) \sim \int f d U$. Projectivity
means that $M_{V}(f)=0$ implies $M_{U}(f)=M_{I}(f)$. 
We have $M_{V}(f)=M_{I}(f * \beta)=M_{I}(f)$ with $\beta \in \mathcal{D}$, 
given $f \in B_{pp}$, that is a sufficient condition for projectivity in the mean, is $f \in B_{pp}$ .

\newtheorem{prop10}[def1]{Proposition (Projectivity in the mean)}
\begin{prop10}
 Assume $U_{1}$ projective over $\mathcal{B}_{pp}(\Omega)$ and that $U_{1} \rightarrow U$ preserves pseudo convexity, then U is projective in the mean,
 that is $M(U f + U^{\bot} f)=M(f)$ over $\mathcal{B}_{pp}(\Omega)$.
\end{prop10}

More precisely, assume $\{ U_{1} f \leq \lambda \} \subset \subset \Omega$ implies $\{ U f \leq \lambda \} \subset \subset \Omega$.
For instance, we can assume $f \in B_{pp}(\Omega)$ and $\mid U_{1} f \mid \leq C \mid U f \mid$.
Given $M(U^{\bot} f)=0$, when $M(U_{1}^{\bot} f)=0$ and when we assume $U_{1}$ projective over 
$\Omega$, through
$f \rightarrow U f$ continuous, we have that $M(Uf)=M(f)$ on compact sets. Conversely, given $M(U f - f)=0$
we have that $M(U_{1}^{\bot} f)=0$ and according to the above $M(U^{\bot} f)=0$.

\subsection{Relative projectivity}

 Assume U analytic with finite Dirichlet integral on W, such that $U^{\bot}=0$ on W and in 
 the same manner for $U_{0},U_{1}$
 harmonic with finite D-integral, so that $U_{1} \leq U \leq U_{0}$ . Thus, we have that
$U_{0}^{\bot} \leq U^{\bot} \leq U_{1}^{\bot}$, that is given projectivity for the both outer sides,
$I-U^{\bot}$ is projective on W, that is \emph{U is relatively projective}. Consider for instance $I-U^{\bot}$ projective on $R(U)^{\bot}$. 

 \newtheorem{def3}[def1]{Definition (Relative projectivity)}
 \begin{def3}
  When a movement $U$ is comparable to projective movements on a subset $W$, we say that $U$
  is relatively projective on $W$.
 \end{def3}

 Assume $U_{1} \subset U$, that is the domain $D(U_{1}) \subset D(U)$. In the same manner assume $D(U_{0}^{\bot}) \subset D(U^{\bot})$.
 Assume further that $U_{0},U_{1}$ are projective according to $D=X_{0} \bigoplus Y_{0}$ and $D=X_{1} \bigoplus Y_{1}$,
 where for instance $D(U_{0})=X_{0}$. Given $X_{0} \cap Y_{0}=\{ 0 \}$ and $X_{1} \cap Y_{1}=\{ 0 \}$, we do not
 necessarily have that $U$ is orthogonal. For instance $X_{0}=X_{1} + Z$ and $Y_{1}=Y_{0} + Z_{0}$. This means that $Z \sim Z_{0}$.

 Assume $U_{j}^{\bot} \rightarrow (-U_{j})^{\bot}$ preserves dimension, with $U_{0} \leq U^{\bot} \leq U_{1}$,
 then $U^{\bot} \rightarrow (-U)^{\bot}$ preserves dimension. Given $f \in \Gamma = \{ U=U^{\bot} \}$, 
 we must have $U \nsim (I - U^{\bot})$.
 Outside $\Gamma$, $U$ does not change character. 

 \subsection{Multipliers}

Concerning multipliers, given $M=1/(X \eta - Y \xi)$, when $-\eta / \xi \sim X^{*} / Y^{*}$, we have that $(X X^{*} + Y Y^{*}) M \sim 1$.
Note that when $-\eta / \xi \sim Y / X$, we have that an infinite D-integral implies $M=0$.
Given $X,Y$ analytic, they can be represented by a Hamiltonian. Assume instead $\log f \in \mathcal{D}_{L^{1}}$, which implies for instance $\phi \sim_{m} \log f$
(\cite{Dahn13}), where $\phi \in H$. $\phi$ can be represented by $e^{\psi}$ with $\{ \psi < \lambda \} \subset \subset \Omega$.
Assume $\widehat{\mu_{\lambda}}(f) (\simeq \mu_{\lambda}(\widehat{f})$), where $M_{\lambda}(f d U)=\int f e^{-i x \dot \lambda} d \mu_{\lambda}$.
One parameter, sequential movements correspond to reduced measures (hypoelliptic d.o). Assume $d \mu$ BV implies
existence of $d v$ reduced and $d \mu=\rho d v$, with $\rho \rightarrow 0$ in $\infty$. 
For instance $\rho=1/Q$, with $Q$ HE polynomial, where we can assume $(d v)^{\bot} \in C^{\infty}$ on $R(Q)^{\bot}$.

Consider $f \frac{\delta \phi}{\delta x}=\frac{\delta f}{\delta x}$, 
we then have $\xi / \eta = (f \xi) / (f \eta)$. That is, given $f \neq 0$, we have that $f X_{V}(\phi)=0$ iff $X_{V}(f)=0$.
Thus, given $f \neq 0$, when the movement is defined by $(\xi,\eta)$, it is sufficient to study the 
movement in phase.
However, a maximum principle for f does not simultaneously imply a maximum principle for $\log f$.

Consider $\Omega=\{ F < \lambda \}$. A \emph{stratifiable domain } (\cite{Dahn13}) is such that 
we have in particular existence of a neighborhood of
$\Omega$, $\Omega'$ such that $\Omega$ is closed in $\Omega'$.
\newtheorem{prop1}[def1]{Lemma (A stratifying multiplier)}
\begin{prop1}
 Assume $\mid U^{\bot} g \mid \sim \mid w_{U} \mid \mid g \mid$, where $1 / w_{U} \in \dot{B}$.
 Then the multiplier $w_{U}$ defines a stratifiable set relative $U$.
\end{prop1}

Assume $\Omega'$ is defined by $\{ w F < \lambda \}$
with $1/w \rightarrow 0$ (in $\dot{B}$), that is $\{ F < \lambda / w \} \subset \{ F < \lambda \}$, close to $\infty$.
Note that completeness for w, is necessary for relative compactness. 
Assume $w(t x,t y)=t^{\sigma}w(x,y)$.
Sufficient for a simply connected continuation (cone continuation ) is $\sigma  > 0$. Assume $\mid U f \mid = \mid w \mid \mid f \mid$,
with $w$ $\sigma-$ homogeneous and $\sigma > 0$, then there is a simply connected continuation. Consider 
$w(x,y) \rightarrow \tilde{w}(u_{1},u_{2})$. When $\tilde{w}$ has cone-continuation, we have that the corresponding
movement  does not change orientation, as $t \rightarrow \infty$.  

\newtheorem{def2}[def1]{Definition (Cone continuation)}
\begin{def2}
 If a movement $U$ can be continued continuously, without changing character, to the infinity. 
 we say that it has a cone continuation.
\end{def2}

Runge's property means that the limit is independent of starting point.
Given a \emph{cone continuation }, we have that dimension is preserved. 
Consider $(x,y) \rightarrow (x,\frac{y}{x}=\rho)$, we then have given $\rho(t x)=t \rho(x)$
a cone continuation . Further $(x,\rho) \rightarrow (x,y)$ defines a convex curve,
that is $y(t x)=t^{2} y(x)$.

\section{Almost orthogonal functionals}
\subsection{Definition}

Transitivity means that $S \sim T$ iff we have existence of a $\sigma \in N(\mathcal{G})$, 
such that $T=S \sigma$ (N is the normalisator).
For instance, assume $V \varphi=\psi$ i $\mathcal{D}_{L^{1}}$
and $\mathcal{F} \varphi=R^{\bot} \psi$ i $L^{1}$, for some $\psi \in \dot{B}$, we then have $R^{\bot} V \sim \mathcal{F}$.
We say that $R,V$ are almost orthogonal, if $R^{\bot} V = R V^{\bot}$. Given $R^{\bot \bot} = R$ and
$R V^{\bot} = V^{\bot} R$, we then have that $R V^{\bot} \bot R^{\bot} V$. In particular, when $R$ is 
the restriction to $\Omega$ in the domain for $d V$ BV, we have that $(d V)^{\bot} \simeq d V^{\bot}$ on $\Omega$.
\subsection{Projectivity}

Given ${}^{t} U^{\bot} \in \widehat{\mathcal{D}_{L^{1}}'}$, we can  assume $R(U)^{\bot}=\{ P=0 \}$.
Consider $E$ as local parametrix to $X_{U}$, where the corresponding U is projective.
Given E is symmetric, then E can be used as orthogonal base. The polar is defined by $\mbox{ ker } E$.
Sufficient for $\mbox{ ker }E=\{ 0 \}$ (modulo $C^{\infty}$), is that $X_{U}$ hypoelliptic in x,y.

Assume $R \gamma = \gamma \mid_{X}$ and $X \cup X_{0}=D$, $F^{\bot}(\gamma)=F(\gamma \mid_{X_{0}})$.
We assume that $F$ has compact support on $X \cap X_{0}$. Let ${}^{t} R F(\gamma)=F^{\bot}(\gamma)$, 
we then have that ${}^{t} R \sim R$ gives a maximal extension. 
Consider $\mathcal{D}_{L^{1}} \subset \dot{B} \subset B$. We can consider $\tilde{\Sigma}$ through
for instance $U_{1}R^{\bot} \varphi \in B$. Note for an oriented foliation $dU=0$, it is necessary  
that $\xi,\eta$ have order 1. 

Consider almost orthogonal movements, on the form $(U \widehat{f})^{\bot} \simeq U^{\bot} f$
that is $\widehat{f} \simeq f^{\bot}$. A closed movement is interpreted as existence of $W \bot I$ 
algebraic, such that $<W \widehat{g},\widehat{g}>=0$, for instance $W \widehat{g}=0$.
Choose $U^{\bot} \in \mathcal{G}$,
that is $d U^{\bot}$ BV and preserves pseudo convexity.
Assume $U$ surjective, we then have ${}^{t}U$ locally  1-1. For instance, when 
$U + V$ projective, we have that ${}^{t}(U + V)$ locally  1-1 and when $V=U^{\bot}$, 
we see that U can not be a spiral. Assume $(Uf)^{\bot}$ a translation domain, that is 
$g \bot Uf$ implies $g=V f$, where $V$ is translation, then the orthogonal is pathwise connected. 
Assume the polar pathwise connected, that is
for every $f,g$ in the polar, we have existence of $U \in \mathcal{G}$, such that $g=Uf$,
that is ${}^{t} U$ locally  1-1. In particular, when ${}^{t} U \simeq U^{\bot}$, then $U^{\bot}$ 
has a representation through
a reduced measure. Projectivity means that the polar can be divided in algebraic components.

\subsection{A separation condition}
 Given $\Gamma$ separates $\gamma$ from $\gamma^{\bot}$, where $\Gamma$ is a simple Jordan curve,
then the domain for $\gamma$ is not simply connected . Given $\Gamma$ is defined by $U \simeq U^{\bot}$,
where $U \rightarrow U^{\bot}$ is defined by a contact transform, we have that  the domain is a simply
connected, if $\Gamma \sim 0$, that is $U^{\bot} \rightarrow 0$ regular continuously and 
the measure for $U^{\bot}=0$ is zero.

Consider $<X(f),\varphi>=
<(\eta_{x} - \xi_{y}) f, \varphi> + <f,{}^{t} X(\varphi)>$. 
In particular, $<f,\widehat{g}>$ and $(\xi,\eta) \bot (-Y,X)$. The condition $(\xi,\eta)$ polynomial , 
means that a movement on f has a corresponding movement in $\widehat{g}$.
Given $(\xi,\eta)$ is defined by G, 
we have when U is harmonic, that $- \Delta G=0$. Further, given $f \in D(U)$ and $d U(f) \bot \varphi$, 
with $\varphi$ analytic, we have that $f \bot R({}^{t} U)$, that is given ${}^{t} U$ preserves 
analyticity and $f \in L^{1}$, then U can be seen as analytic.

\section{Concepts from spectral theory}
\subsection{Numerical range}
Assume \emph{ the complement to NR}
is generated by one single movement, we then have $\Gamma \subset NR$. 
Assume $<U f, V f>=<{}^{t} V U f,f>$ defines numerical range NR. Given ${}^{t} V U$ a bounded
operator, we have that $\sigma({}^{t} V U) \subset NR$.  Consider $d {}^{t} VU=\frac{d {}^{t} V}{d U} (\xi d x + \eta d y)$ and 
$(\xi,\eta) \rightarrow (\xi_{x} + \eta_{y},\eta_{x} - \xi_{y}$),
that is harmonic conjugates. 
Note that when $U \rightarrow {}^{t} U$, such that $\xi_{x} + \eta_{y}=0$, 
the character of movement is preserved.
Note that when $U$ is not ac, 
there are examples of U monotonous with $d U=0$, but 
$U \neq \lambda I$.
Define $I(\sigma)$ as $f$, where the movement changes character. 
Note that given $d U = \alpha d U_{1}$, then $\alpha$ is regular
outside $I(\sigma)$.
\subsection{Spectral resolution}
Assume $\exists U \in \mathcal{G}$, such that $V f=f(U^{\bot} \gamma)=0$. We then have given 1/f a.p., 
that Vf defines invariant sets. A very regular boundary means existence of a regular approximation, 
that is we allow spiral approximations.
Given the inverse to $U^{\bot}$ continuous and ker f connected , the invariant set is connected .

Given $V f \in C^{\infty}$, we have  $\frac{1}{R} \log \mid \widehat{Vf} \mid < 0$, when R large.
For instance, $V f = \int \int_{\rho < \lambda}f e^{i(x-y) \dot \xi} d \xi d y$, for $d V^{\bot}=\rho d I$, 
that is we can 
represent V as regularizing with kernel d V in $C^{\infty}$ outside $\Gamma$. Necessary for this representation, 
is that $I \prec \prec V^{\bot}$, that is we associate $(d V)^{\bot}$ to a reduced measure (\cite{Nilsson72}).

Assume further $V=U_{1}$ harmonic at the boundary, we then have (point wise topology) that $V(fg)=V(f)V(g)$ 
and $\log V f \sim V \log f$. Consider the condition $\{ U,V \}=\xi \eta' - \xi' \eta=0$, that is 
$\xi / \eta \sim \xi' / \eta'$. We then have $\frac{d V}{d U} \sim \frac{\xi' \xi + \eta' \eta + i(\xi \eta' - \xi' \eta)}{\mid \xi + i \eta \mid^{2}}$,
that is given $\xi' / \eta' \rightarrow -\eta / \xi$ projective or $\{ U,V \}=\{ U,V^{\diamondsuit} \}=0$,
we have that $d V \bot d U$. Given $d (U-I)=0$ implies $U=U_{j} \in \mathcal{G}$, $j=1,\ldots, k$, 
defines the boundary (subset of polar), we have that when the boundary has order 0, 
then $k=1$. Thus, we have that $V$ is
harmonic at the boundary and that $d U \bot d V$ there.

A Banach space is Hilbert iff $\parallel U f+ U^{\bot} f \parallel^{2} + \parallel U f-U^{\bot} f \parallel^{2}=
\parallel U f \parallel^{2} + \parallel U^{\bot} f \parallel^{2}$,
$\forall U f, U^{\bot} f \in B$. Assume $U,U^{\bot}$ are symmetric, normal operators, with $U U^{\bot}=U^{\bot} U$.
We then have $U + i U^{\bot}$ is normal
iff $U U^{\bot}=U^{\bot} U$. Further, $U \bot U^{\bot}$ if $U U^{\bot}=U^{\bot} U=0$. Assume 
$U-I \in C^{\infty}$, with $\mbox{ ker }U=\{ 0 \}$.
Further $\parallel (U + I) f \parallel = \parallel (U-I) f \parallel$. Given $(U-I) f=0$, we have that
$f=Uf=0$, that is a ``symplectic'' model. Assume $d U \bot d V$ with $d U$ BV. Necessary for a 
classical spectral resolution, is $d V$ BV. The condition $I \prec V$ (sufficient is that V is algebraic), 
does not imply that U is projective.

\subsection{Orthogonal base}

Dirichlet problem is to determine $u \in C^{0}(\overline{\Omega})$, with $\Delta u=0$ on $\Omega$
and $u=f$ on $bd \Omega$ (unique through the max-principle). Assume $\frac{\delta u}{\delta x}=\frac{\delta v}{\delta y}$
and $\frac{\delta u}{\delta y}=-\frac{\delta v}{\delta x}$.  Then the condition is that $\frac{\delta^{2} v}{\delta x \delta y}=
\frac{\delta^{2} v}{\delta y \delta x}$. Simultaneously, we have that $\{ u,v \}=(\frac{\delta u}{\delta x})^{2} + (\frac{\delta u}{\delta y})^{2}$.
The Dirichlet problem for movements is existence of U harmonic on $\Omega$ with $U=U_{j}$ on $\Gamma$
that is $d U = \beta d U_{j}$ with $\beta = const$ at the boundary $\Gamma$. Thus, there is a Hamiltonian G
symmetric at the boundary, corresponding to V. Given a maximum principle and U continuous at $\Gamma$, 
U can be determined as $U_{j}$ uniquely. 
Given that U is harmonic, we have  an orthogonal base iff V is symmetric and projective. Presence of 
a max-principle, is dependent on if
the domain is contractible and does not simultaneously imply a max-principle in phase. 
Given $\mid U e^{\phi} \mid \leq e^{\mid w \mid \mid \phi \mid}$, with $\mid w \mid$ finite, gives a
upper limit. Presence of a max-principle means that the limit is reached. In particular
$\mid d U(f) \mid \leq \mid X(f) \mid \leq \max(\mid \xi \mid,\mid \eta \mid) \mid d f \mid$.
Consider $d U=0$, given $U$ ac, then $U$ has a max-value that is reached. Or we can assume $d U=\beta d U_{1}$
with $d U_{1}$ harmonic (min and max are reached) with $\beta=const$ close to the boundary. 
Concerning \emph{flux}, Assume $U^{\bot} \rightarrow U^{\diamondsuit}$ continuous
where the movement is related to axes, that is $L^{\bot}=\rho L^{\diamondsuit}$ with $\rho=const$, 
we then have that 
$\int_{\Gamma} d U^{\diamondsuit}=0$ implies $\int_{\Gamma} d U^{\bot}=0$. In particular when 
$d U^{\bot}=\beta d U^{\diamondsuit}$ with $\beta=const$ close to the boundary, we have that  
flux is preserved. In the case when U is defined in the phase, we assume  $U \phi \in L^{1}$,
that can be approximated by $U \phi \in H$. In particular if $d U \phi,d U^{\diamondsuit} \phi$ harmonic, 
is $d U \phi$ analytic. 
\section{Boundary conditions}
 
 Assume the polar $C$ is defined as zeros to a holomorphic function. Given $\tilde{C}=\{ (x,y) \quad U^{\bot}=0 \}$
 and $\psi : \tilde{C} \rightarrow C$ continuous and proper, if $U^{\bot}$ analytic, we have that 
 $\tilde{C}$ is removable iff we have a global base for the corresponding ideal. In particular 
 when $U^{\bot}$ algebraic, we have that $\tilde{C}$ is removable. Note that given the boundary 
 defines a strictly pseudo convex set according to the above,
 the normal can be given locally by polynomials  (\cite{Oka60}). 
 
 Assume M the cylinder web and that dim = $dim_{\mathcal{G}}$, that is rank is taken relative the 
 group of movements. For instance $U_{S} \notin \mathcal{G}=\mathcal{G}_{1}$, we have that $dim_{\mathcal{G}} M=0$,
 but M is two dimensional.

Assume $U$ projective and analytic with $\mid \xi \mid^{2} + \mid \eta \mid^{2}$ finite on $E = \complement \Omega$, we then have 
that $d U=0$ on E implies $U^{\blacktriangle}=0$. When $\Omega$ planar $\in \mathcal{O}_{AD}$ (\cite{AhlforsSario60}) we have that
every single valued and linear u is constant and conversely, that is $U^{\blacktriangle}=0$ defines E. When E is closed and $u \in AD(nbhd E)$, we have that  
E is removable iff $\Omega \in \mathcal{O}_{AD}$. 
A global model has a removable polar, this is assumed invariant for algebraic changes of local 
coordinates.

 \subsection{Co dimension 1}
 Assume $\Gamma$ points where the movement changes character, we then have $U_{S}=U^{\bot}_{S}$ 
 is included in $\Gamma$. 
Assume $d U^{\bot}=\alpha d U$ and $\tilde{\Gamma} = \{ \alpha=const \}$, Thus, the boundary $\tilde{\Gamma}$
is not oriented with respect to $x,y$, when $\tilde{\Gamma}$ contains a spiral. The co dimension for V and $\tilde{\Gamma}$ 
is the maximal dimension for $f \in C^{\infty}$ 
that approximates $\tilde{\Gamma}$, such that $V f \in \mathcal{D}_{L^{1}}$. 
When $d V$ does not change character, 
the complement to $\tilde{\Gamma}$ must be connected, 
that is $0 \neq f,g \in \tilde{\Gamma}^{\bot}$ implies $g=V f$. When $V$ is not reduced, $Vf \in L^{1}$ 
does not imply $f \in L^{1}$. Using monotropy, we can  define the order for $\tilde{\Gamma}$ as the number
of defining (linearly independent) functions. When $d V$ BV, we can assume $\tilde{\Gamma}$ has co dimension 1.

 \subsection{Reflexivity}
 When $\mathcal{G} \times \mathcal{G}^{\bot}$ is assumed very regular, we assume  existence of regular 
 approximations, for instance $d U_{1}$ with $U_{1} \neq U_{1}^{\bot}$. Simultaneously, we can  have 
 $\exists j \quad d U_{j}=0$ and $U_{j} \simeq U_{j}^{\bot}$. 
 Consider $\mathcal{G}^{\bot}$ as the completion of $\mathcal{G}^{*}$ to $L^{1}$.
 Thus $\exists U \in \mathcal{G}$,
 such that $U \bot \mathcal{G}^{\bot}$. Assume $U^{\bot}=(I-U) + V$ is closed,
 we then have $U^{\bot \bot} = (I-U)^{\bot} + V^{\bot}=U - V + V^{\bot}$, that is $U^{\bot \bot} - U = V^{\bot} - V$.

\newtheorem{prop4}[def1]{Proposition (Diagonal continuation)}
\begin{prop4}
 Assume $U$ a reflexive, but not projective movement. Then the continuation $U + V$, is diagonal,
 that is $V^{\bot}=V$.
\end{prop4}

When -V=I, we have that $U^{\bot}=-U$. Given $V=V^{\bot}$ implies V=0, 
we have  that reflexivity implies projectivity. Thus $V=V^{\bot} \neq 0$ is polar for reflexivity. 
Existence of $dV \neq 0$ implies existence of non-trivial boundary $\tilde{\Gamma}$, that is polar for 
projectivity.

Assume $d U_{i}=\alpha_{ij} d U_{j}$. Given $\alpha_{ij} \rightarrow 0$ in $\infty$, we can define
$U_{i} \prec \prec U_{j}$, that is inclusion for the corresponding space of integrable functions, 
for the limes to be an isolate point, the $\alpha_{ij}$ must be regular, that is
the concept of co dimension is dependent of the inclusion condition.  
Leaves are connected components, given $U_{i}$ continuous. The leaves are assumed defined by 
movements of the same character and of constant rank. 

 \subsection{Multipliers}
 Assume $\frac{\frac{1}{xy}}{V}=\frac{x}{y} + \frac{y}{x}$, where $1/V=x^{2} + y^{2}$. Given $\rho=\frac{1}{xy}$, we have that
$\rho \frac{1}{M}=\frac{1}{N}$, that is given N exact then M is exact (\cite{Lie91}). Further $1/\rho(x,y)=\rho(\frac{1}{x},\frac{1}{y})$. 
Given $\rho^{*} \widehat{N}=\widehat{M}$ where $\rho^{*}=\widehat{H}$, where H is Heaviside, we have that
$M=H* N$, that is M can be represented with support in a rectangle.
Let $R_{1}(x,y)=(y,x)$ and
$R_{2}(x,y)=(x,-y)$, we then have $\diamondsuit=R_{1}R_{2}$. Note that if $\rho \rightarrow \rho^{*}=\frac{1}{x^{*} y^{*}}$,
we have that $\frac{d^{2} H}{d x d y} \widehat{*} F= \delta_{0} \widehat{*} F$. Note that $\int_{\Omega} 1 d x^{*} d y^{*}$ finite
assumes $\Omega \subset \subset \mathbf{R}^{n} \times \mathbf{R}^{n}$. 

Consider $<R T,\varphi>=<T,{}^{t} R \varphi>$, we then have if R restriction, then ${}^{t} R$ is extension
and conversely. Assume $f \in (I)$ implies 
$f=f_{0} + f_{0}^{\bot}$, we then have $R f_{0}=f_{0}$ iff $R^{\bot} f_{0}^{\bot}=f_{0}^{\bot}$.
Further $(R-I) f_{0}=0$ iff $(R^{\bot} - I) f_{0}^{\bot}=0$. Further $R^{\bot} f_{0}=0$ iff $R f_{0}^{\bot}=0$
that is $R^{\bot} f = f_{0}$ and $(R + R^{\bot}) f=f$. For the spiral, we have that $U f= U^{\bot} f$
with $U^{\bot} f_{0}= U f_{0}^{\bot}$, that is given $U$ projective over $f_{0}$, then $U^{\bot}$ is projective
over $f_{0}^{\bot}$.

\section{Movements modulo $C^{\infty}$}

\subsection{Convexity}

When $X dt=dU$ has positive derivative ( one parameter), that is a strictly convex measure, then dU is monotonous 
and has a division in an ac and a singular part. 
 Consider $f=\beta e^{v}$ and a set such that $v =0$ and $\log \beta \leq \epsilon$.  
 Assume $\beta \neq 0$. Consider $f^{2} \sim \beta^{2} e^{v}$ that is if $f \sim e^{v}$,
 then f can be seen as geometrically radical. Consider $(f_{1} + i f_{2})^{2} \rightarrow \infty$, 
 a sufficient condition is that $1 / ( \frac{f_{1}}{f_{2}} + \frac{f_{2}}{f_{1}}) \rightarrow 0$, 
 that is $f_{1} \prec \prec f_{2}$ or $f_{2} \prec \prec f_{1}$, that is one sidedness. 
 Spirals approximate multivalent surfaces, 
$(U_{S} - U) \rightarrow 0$.
Note that $U_{S}$ can be seen as monotonous i $(s,t)$,
which does not imply monotonous in $s,t$. 
 
 Assume the movement is considered modulo $C^{\infty}$, that is $U^{\bot} f =0$ implies $(I-U)f \in C^{\infty}$.
 Note that if V is projective, $(V - I)^{2} \simeq (V - V I) + (I - V I)$ and in $H$, we have that
  $\simeq (I -V)$, that is \emph{$(V - I)^{2} f \in C^{\infty}$ iff $(V - I) f \in C^{\infty}$}. 
  
  \newtheorem{prop8}[def1]{Proposition (The orthogonal is radical)}
  \begin{prop8}
   Assume $V^{\bot} f \subset C^{\infty}$ and that V is projective over $f \in H \cap \mathcal{D}_{L^{1}}$.
   Then, $(V - I)^{2} f \in C^{\infty}$ iff $(V - I) f \in C^{\infty}$.
  \end{prop8}
 
 $\Sigma X_{j}^{2}$ is convex. Assume $d V \simeq d U^{2}$ and $\Sigma X_{j}^{2} \simeq X_{V}$.
Thus,  given $d V$ convex we have that $d V(g)=0$ implies $(I-V) g =0$ on compact sets. 
Note $(U+I)(U-I)=U^{2} - I$ (cf. defect indexes). However we do not have
d V convex implies V projective.
Assume $C=\omega(U)$ points where the movement changes character. Obviously, 
$\omega(U) \neq \omega(U^{2})$. Given $g$ symmetric around 0, we have that $(U-I)g=0$ iff 
$(U+I)g=0$, that is $\omega(U) \simeq \omega(U^{2})$.
Note that ac is not a radical property, that is $d U^{2}$ ac on compact sets, does not imply 
that $d U$ is ac on compact sets.
Assume $d U=0$ implies $U-\rho(t) I=0$, with $\frac{d \rho}{d t} \neq 0$, for instance with $U - \rho I \leq \epsilon$.
Further, if $(U^{2}-I)g=0$ in isolated points, then $(U -I)g=0$ can still contain a segment
of functions (a positive dimension for the space of eigen vectors).  Assume $g \bot d U$ and 
$\widehat{g} \bot d U^{\bot}$. 
Given $U,U^{\bot}$ ac and $U^{\bot \bot} \simeq (I-U^{\bot})$ with $U$ reflexive, we get 
$\parallel U^{\bot} g \parallel = \parallel U \widehat{g} \parallel$.
If further $U^{\bot}$ is normal, we have that $\parallel U g \parallel = \parallel U \widehat{g} \parallel$.
Given w complete, it is sufficient to put the condition
for $\widehat{U \frac{\delta f}{\delta x}} \simeq w x \widehat{f}$, where 
$x w \sim w'$ and  $1/w' \rightarrow 0$ in $\infty$, that is the condition for $I \prec \prec X(\log f)$.

\subsection{Projectivity modulo $C^{\infty}$}
Assume $(d U)^{\bot}=d V \subset C^{\infty}$ (point wise topology). Given $V(\varphi)=v*\varphi$, we have that
$\frac{d}{d x}V(\varphi)=v * \frac{d \varphi}{d x}=V(\frac{d \varphi}{d x})$. Assume $(J)=\{ \varphi \quad \widehat{\varphi} \in \mathcal{D}_{L^{1}} \}$
since $\widehat{\delta_{0}}=1$, we have that over $J$ that $\delta_{0}$ is algebraic. 
Given $U f \in (J)$, we have that  UI=IU over (J). Assume $U$ is very regular, that is $U f=f - C^{\infty}$, 
we then have that $R(U)$ is dense in $\mathcal{D}_{L^{1}}$, that is given $f \in \mathcal{D}_{L^{1}}$, $U f - f \rightarrow 0$ i $\mathcal{D}_{L^{1}}$,
as $U \rightarrow I$. Assume $U \in \mathcal{G}_{HE}$ with $UI=IU$ implies $(I-U) \in C^{\infty}$, then 
$U(I-U)=(I-U) =0$ modulo $C^{\infty}$, that is $U^{2}=U$ modulo $C^{\infty}$.

Given U, we can determine a domain, where U is projective and $U^{\bot} f \in C^{\infty}$, $\forall f$. 
Assume  $U^{\bot} f = g \in L^{1}$,
we then have outside the polar that $g \sim U_{1} f$ and $d U = \alpha d U_{1}$, $d U^{\bot}=\beta d U_{1}$.
Starting from 
graph norm, there is through Riesz-Thorin a maximal domain for projectivity (\cite{Dahn19}). 

Given $\gamma$ regularizing in $L^{1}$, we have that $\{ \gamma < \lambda \} \subset \subset \Omega$.
Note that over $\Sigma_{V,\lambda}=\{ M_{V}(f) < \lambda \} \subset \subset \Omega$, 
$\int_{\Sigma_{V,\lambda}} f d V < \lambda \int_{\Sigma} d x$.
Given $f_{0}$ very regular, we have that $\{ P(\xi) f_{0} < \lambda \}$  $ \subset \{ P(\xi) < \lambda \}$,
that is assume $Q$ hypoelliptic, with $f_{0} \sim 1/Q$. We then have $\mid \xi \mid^{\sigma} \leq \mid Q \mid$
and $0 < \sigma \leq 1$, we have that $\{ P / \mid \xi \mid^{\sigma} \leq \lambda \} \subset \{ P f_{0} \leq \lambda \}$
$\subset \{ P / \mid \xi \mid^{m} \leq \lambda \}$, where $m$ is the degree for $Q$.

\emph{Assume $U + U^{\bot}=A$ with $A - \delta \in C^{\infty}$}. Given $U \neq U^{\bot}$, we have that $A (U,U^{\bot}) - \delta \in C^{\infty}$.
Assume $U$ projective in mean. Given $U \bot V$ (with $U^{\bot} \sim V$), we have that U is 1-1 in graph norm. 
Assume $U f \in \dot{B}$ and $<U f,V f>=0$, we then have given $f \in \dot{B}$, ${}^{t} V U \in \mathcal{D}_{L^{1}}'$. 

\newtheorem{lemma14}[def1]{Lemma (Reduced movements)}
\begin{lemma14}
 Consider $U + U^{\bot}=A$, with $A - \delta \in C^{\infty}$. Then U is projective (modulo $C^{\infty}$)
 outside the kernel to $A$.
\end{lemma14}
Note that outside the kernel, $A$ is reduced modulo $C^{\infty}$ and ker A $ \supset \Gamma$. 
Given $U$ projective modulo $C^{\infty}$ outside the kernel to A, that is $U + U^{\bot} - I \subset C^{\infty}$
and given $U^{\bot} f \subset C^{\infty}$, we have that $U - \lambda_{j} I \rightarrow U - \lambda I$
as $0 < \lambda_{j} \rightarrow \lambda$.
We have that, $A=A(U,(1-U))$ 
and when $A=0$, we have that $U = - U^{\bot}$, that is a spiral domain . Assume $U$ projective on $X$, that is $A=I$ on $X$.
Assume A very regular on $X^{\bot}=X_{0}$. we then have $A=A(U)$ on X and $A=A(U,U^{\bot})$ on $X_{0}$.
Thus, the co dimension is dependent on projectivity for U.
\subsection{Continuations}
Assume instead of projectivity, $U^{2} \varphi - U \varphi = -U (I-U) \varphi \in C^{\infty}$. 
For differential operators, when solutions are considered modulo $C^{\infty}$, it is sufficient 
to consider operators of real type.
Let $T_{\beta}(\varphi)=\int \varphi * \beta = T(\varphi * \beta)$ analytic, where $\varphi * \beta \in \mathcal{D}_{L^{1}}$
analytic and $\alpha,\beta \in \mathcal{D}$. 
Consider $U \sim P(D) f_{0} \in \mathcal{D}_{L^{1}}'$, given that
$f_{0}$ very regular, we have that $T_{\alpha}(U \varphi)-U T_{\alpha}(\varphi)=(U \varphi) * \alpha - \varphi * (U \alpha) \in C^{\infty}$. 
Given $U + U^{\bot}=I$ over X, the equation can be continued to $U + U^{\bot}=I + W$, with 
$W \in C^{\infty}$ over $X_{0}$. \emph{When $X \bigoplus X_{0}=L^{1}$, Consider $X=\mathcal{D}_{L^{1}}$, we then have $\widehat{X_{0}} \subset \dot{B}$}.
Assume that 
$X_{0}$ is defined by $\frac{d^{j}}{d x^{j}} \varphi \in L^{1}$ implies $\frac{d^{j}}{d x^{* j}} \widehat{\varphi} \in \dot{B}$, 
$\forall j$, we then have $ \widehat{U \varphi} \in \dot{B}$. Assume $U^{\bot} \widehat{\varphi}=\widehat{U \varphi}$.
When $\varphi \in \dot{B}$, then U can be represented in $\mathcal{D}_{L^{1}}'$, that is $\sim P(D) f_{0}$, 
where $f_{0}$ very regular. If we assume $P$ partially hypoelliptic, then the kernel to $f_{0}$ can be 
seen as very regular.
 Note $\int \mid d U \mid$ finite
implies finite D-integral. Given $U^{\bot} \varphi \in \dot{B}$
we have that $\int \mid d U^{\bot} \mid$ locally  finite. 

\newtheorem{lemma15}[def1]{Proposition (Implicit projectivity)}
\begin{lemma15}
 Assume $U$ is defined by $U^{\bot}$, that is over $\mathcal{D}_{L^{1}}$ $\widehat{U f}=U^{\bot} \widehat{f}$.
 Then, $U$ is projective over $\mathcal{D}_{L^{1}}$ as long as $U^{\bot}$ is projective over $\dot{B}$.
\end{lemma15}
Sufficient is to note that $\mathcal{D}_{L^{1}} \subset \dot{B}$.
Assume $U^{\bot}$ projective over $\dot{B}$ and that U is defined through $U + U^{\bot}=I$ in 
$\dot{B}$. Since $\mathcal{D}_{L^{1}} \subset \dot{B}$, then U must be projective in $\mathcal{D}_{L^{1}}$.
For projectivity in $L^{1}$, we must assume for instance $(I-U) \varphi \in \dot{B}$.

Note that we assume $U F(\gamma)(\zeta)=F({}^{t} U \gamma)(\zeta)=F(\gamma)(\zeta_{T})$.
Given $\zeta$ reduced in $x$, we have  $\mid x \mid^{1/N} \mid 1 / \zeta \mid \leq C$, that is $x$ 
of algebraic growth in $\zeta$. When $y=\varphi(x)$, then $\varphi$ can be
chosen as projective, that is a standard complexified situation. 
Under the condition above, given $U^{\bot} \gamma \subset C^{\infty}$ and $U^{\bot \bot} \gamma \subset C^{\infty}$
it is sufficient to consider reflexive one parameter movements, to define $\mathcal{G} \times \mathcal{G}^{\bot}$.

Given U analytic,
then $d U=0$ is a closed form and when $d U = 0$ is reduced, the form is exact. Given $I \prec U$ then U preserves compact sub level sets locally
and given $I \prec \prec U$, U is locally reduced. When U is reduced, we have existence of $U^{-1}$. 
Given $\{ U,V \}=0$ with $\xi,\eta \neq 0$, we have that V can be represented as an invertible distribution. 
Given $U^{2}=I$ then U is involutive and $U \simeq U^{-1}$. Given $U^{2} =U$, then U is projective. 
Given $\varphi \in \dot{B}$ implies $U \varphi \in \dot{B}$, we have that $\exists U^{-1} \in \mathcal{D}_{L^{1}}'$. 
When $U^{2}=I$, we have that
$\{ U,U^{-1} \}=0$. Assume $UV=I$ we then have $\frac{\delta}{\delta x} UV=V \xi_{U} + U \xi_{V}=\frac{\delta}{\delta x}I$.
Further $X_{V}(\xi_{U}) + X_{U}(\xi_{V})=X_{I}$

A condition necessary  for hypoelliptic, is that the derivatives are strictly weaker, 
for instance $\xi,\eta \prec \prec U$. Further when the movement is considered in
the phase, given $\xi$ reduced, we have that $e^{1/\xi} - 1 \rightarrow 0$ i $\infty$. When $\xi$ 
linear in $x,y$, we have that $e^{\xi} \rightarrow \lambda I$, as $x,y \rightarrow 0$. 
Given $V \prec \prec U$ and $W \bot U$, then $V \prec \prec W$, that is
consider $(I_{U})^{\bot} \subset (I_{V})$. Given $\xi_{1}/\xi \rightarrow 0$,$\eta_{1} / \eta \rightarrow 0$ in $\infty$,
we have  $d U_{1} \prec \prec d U$ and so on. 

\subsection{Localization's}

Assume $\Phi_{1}(f)=\xi \frac{\delta f}{\delta x} + \eta \frac{\delta f}{\delta y}$ and assume locally 
$E_{1} \Phi_{1} f = \Phi_{1} E_{1} f = f$. Given $E_{1} (I) \rightarrow (I)_{ac}$
surjective, we have that $\forall g \in (I)_{ac}$, we have existence of $f \in (I)$, such that $E_{1}(f)=g$.  
Assume $\{ E_{1},E_{2} \} =0$, where the derivatives are taken in the weak sense. Assume $E_{2}(\frac{\delta f}{\delta y})=
E_{1}(\frac{\delta g}{\delta y})$. We then have $\{ E_{1},E_{2} \}(f)=$ $E_{1}(\frac{\delta f}{\delta x})
E_{1}(\frac{\delta g}{\delta y}) - E_{1}(\frac{\delta f}{\delta y})E_{1}(\frac{\delta g}{\delta x})$.
Assume $E_{1}(fg) \sim E_{1}(f) E_{1}(g)$, we then have $E_{1}(\{ f,g \})=0$. Note that when $E_{1}$ has real type, 
we have that we can choose $fg$ real.
When $fg$ polynomial , we have that $E_{1}(fg)=0$ implies the kernel to $E_{1}$ is $\equiv 0$ 
on a connected  set.
When $f,g \in C^{\infty}$, we have that
$I(fg)=I(f)I(g)$. When $E_{1}=\delta + \gamma$ very regular with trivial kernel, we can  choose $\gamma$, 
such that $f \gamma_{f} + g \gamma_{f}=0$ and $\gamma_{fg} - \gamma_{f} \gamma_{g} \in C^{\infty}$. 
Thus we have that modulo $C^{\infty}$, that $E_{1}(fg) - E_{1}(f) E_{1}(g)=0$. Given $E_{1}$ has non-trivial
kernel, wee must assume the kernel algebraic, that is $E_{1}(fg)=0$ iff $E_{1}(f)=0$ or $E_{1}(g)=0$.

Assume $E_{1}(f)=\widehat{T}(f)$ and consider $\int T(f-g) T(g) d g$. Given $f \sim g$ (conjugation)
we have  $=T(0) \int T(g) d g=T(0) \widehat{T}(0)$, that is given $f \sim g$ with $f,g \in R(E_{1})$ 
is the conjugation corresponding to $\Phi_{1},\Phi_{2}$, we have that $E_{1}(f) E_{1}(g) \sim T(0) \widehat{T}(0)$, 
that is $E_{2}(f) \sim E_{1}(g)$ implies $E_{1} \sim E_{2}$ according to $T \rightarrow \widehat{T}$. 
The equation $\Phi_{1} E_{1} f \sim f$ is interpreted so that f is a symbol, $E_{1}$ is localizer and 
$\Phi_{1} f=0$ defines a regular approximation, $E_{1} ({}^{t} \Phi_{1} f) \sim E_{1} (U^{-1} \overline{\delta} U f)$,
such that the condition ${}^{t} \Phi_{1} f=0$ implies $U f$ analytic.  

Note $\Phi_{1}$ hypoelliptic does not imply ${}^{t} \Phi_{1}$ preserves hypoellipticity (spiral). 
But given $E_{1} f - f \in C^{\infty}$, we have that  ${}^{t} \Phi_{1} E_{1} f-f \in C^{\infty}$. 
Note that $\Phi \rightarrow {}^{t} \Phi$ is dependent of algebraicity for $\xi,\eta$. When 
$E_{1} \Phi_{1} fg=fg$, we have that $E_{1}(g \Phi_{1}(f) + f \Phi_{1}(g))$ $=<E_{1}, g \Phi_{1}(f) + f \Phi_{1}(g)>=<I,fg>$, 
that is given $E_{1}$ linear, we can  write $<E_{1}, \frac{\Phi_{1}(f)}{f} + \frac{\Phi_{1}(g)}{g}>=I=Exp(0)$.
Assume $\widehat{E_{1}}(\phi)=E_{1}(f)$, where $f=e^{\phi}$. Given $\frac{\Phi_{1}(f)}{f}=\Phi_{1}(\log f)$,
we have that $<\widehat{E_{1}}, \Phi_{1} \log f + \Phi_{1} \log g>=0$. Given $<E_{1},E_{2}>=0$, we have that 
$E_{1} \bot E_{2}$ implies $E_{1}^{\bot} \simeq E_{2}$. Given $E_{2}$ projective, we have that $E_{2}^{\bot} \sim I - E_{2}$.

\section{Linear independence}

\subsection{Desingularization}
Note the following problem:  $U_{1}^{\bot} \gamma=U_{2}^{\bot} \gamma$
implies $\gamma=0$. This depends both on projectivity for $U_{j}$, and on the envelop property, 
that is the property that for any movement $U \in \mathcal{G}$,  $U=V^{\bot}$ for some movement 
$V \in \mathcal{G}$. 
Assume that the diagonal $U=U^{\bot}$ is excluded.
and that movements are reflexive, we then have $ V \neq U$ on the diagonal.

Given three equations $X_{j}d y - Y_{j} d x=0$, $j=1,2,3$, with integrals u,v,w and $\lambda u + \mu v + \nu w=0$, where $\lambda,\mu,\nu$
are non zero constants. \emph{Consider the equations $d u=0$,$d v=0$ and $d v - \Phi(u,v)d u=0$. 
Since $\frac{\delta \Phi}{\delta u} + \frac{\delta \Phi}{\delta v}=0$
we have $\Phi=\Phi(u-v)$} (\cite{Lie91}  Ch. 5, Sats 7). Assume now $d U(f)=\{ f,G \}$, for some function 
$G$. We thus have $A(G)=X \frac{\delta G}{\delta x} -Y \frac{\delta G}{\delta y} = d U(f)$, 
given f a Hamiltonian.

Assume $dV=\Phi(U,V) dU$ implies $\Phi(U,V)=\Phi(U-V)$, for $V=U^{\bot}$. When $\Phi$ is a polynomial in 
$U,V$ with $\Phi(I)=0$, we have that
$\Phi(U-V)=0$ implies $U-I \sim V$ locally. 
Thus, $\int d V=\int \frac{d V}{d U} d U=\int \Phi(U-V) d U$. 
When $U$ increasing implies $V$ decreasing, we have that $\frac{d V}{d U} < 0$ implies $\{ U \leq V \}$ compact.
Assume $W=U^{\bot} - U$ and that $\{ U \leq U^{\bot} \}=\{ W \geq 0 \}$ is compact. Given $R(W) \subset C^{\infty}$,
we can  assume $W(fg) \sim_{0} W(f)W(g)$, that is $W$ is essentially algebraic.

\subsection{Factorization}

Consider systems where the radical is isotermal.
Consider $d U = X(f) d t=\frac{1}{M} d t$. For instance $Y/X \sim - \xi / \eta$. When $X,Y$ are real, we can compare 
with the condition of finite D-integral, that is normal surfaces. Consider a system where $(X^{v},Y^{v})$ real 
and with finite D-integral  for some v,
that is the radical to the system defines a normal surface. The corresponding $M^{1/v}$ is integrable, when M is integrable.

Assume $T^{\bot}=U$ analytic on $\Xi$ and $T$ analytic on $\Omega$. Consider $\Xi \cup \Omega$.
Let $T - U = S$ on a domain 
$\Pi$ where $S=0$, that is $T=U$ over $\Pi$, that is $T^{\bot} \sim T$ over $\Pi$, for instance a 
spiral. Given $g=Q f$, for some
analytic movement  Q and $T f + U g=S f$, we assume  Q (interpolation) preserves analyticity $\Xi \rightarrow \Xi$.
Thus,  we have existence of a continuous (ramifier), such that $Q f = g \rightarrow f$, as $Q \rightarrow I$.

\subsection{Multipliers}

Concerning the concept of dimension, note that through the condition on very regular boundary, we have existence of 
$U_{j}$ analytic over some domain $\Omega_{j}$. Conversely, assume every subset $\subset \Omega$ is 
a domain  for analyticity for some movement $U_{k}$ (or several), that is the domain can be generated 
by an analytic movement in $\mathcal{G}$. 

Assume $1/M \sim X \eta - \xi Y$. Given $- \eta / \xi \sim Y^{*} / X^{*}$, we have that $M=0$ iff $D(f)=\infty$.
Given $M = \rho N$ with $\rho$ regular, we have simultaneously $\frac{1}{N}= \rho \frac{1}{M}$. Simultaneously, given $M=\infty$ implies $D=0$.
Consider $X_{U_{S}}$ factorized over $\Phi=(X,Y) \rightarrow (X^{*},Y^{*}) = \Phi^{*}$,
where $\Phi,\Phi^{*}$ are assumed in involution. Note that to determine $X_{U_{S}}=0$, it is sufficient 
to consider $\Phi=const, \Phi^{*}=const$. The spiral closure that we are considering is given by
$X_{U_{S}}=const$. Consider $d \Psi : d U \rightarrow d U^{\bot}$. Given absolute continuity, we have that
$d U=0$ is mapped on $d U^{\blacktriangle}=0$, that is closed forms are mapped on to closed forms.
Given $d^{2} U=0$ implies $d U = c d t$ where the right hand side does not contribute to the integral, 
the movement is densely defined. The condition $d^{2} U = 0 \rightarrow d^{2} U^{\bot}=0$ can be compared 
with $\delta X=0$, that is $U$ harmonic implies $U^{\bot}$ harmonic.

 \section{Symmetry}
Concerning the maximum principle; $\delta_{0}$ reflects the support through 0, that is $TI=IT$ implies that $T$ has support
in a ball formally. Assume $R_{1}$ reflection through the x-axle and $R_{2}$ reflection through the y-axle
(\cite{Dahn14} the two mirror model). We then have $\delta_{0}=R_{1}R_{2}$. Given $TR_{j}=R_{j}T$ $j=1,2$
and $R_{1}R_{2}=R_{2}R_{1}$, then T has support in a ball.
Assume ${}^{t} R_{1}=R_{2}$ and $T R_{1}=R_{2} T$, we then have $T R_{1} R_{2}=R_{1} T R_{1}=R_{1} R_{2} T$,
that is an algebraic relation. 
Given scaling orientation, note
that symmetry as $B \rightarrow 0$ does not imply symmetry $B \rightarrow \infty$. 
Note that when $T(x,x^{*})$ symmetric, with $x^{*} \simeq 1/x$, we have that $T$ is symmetric as $B \rightarrow 0$
and as $B \rightarrow \infty$. 

Consider $d F - d L=0$, where $L$ linear. Given $(F-L)$ ac, we have that $F$ is linear.
Assume $d y / d x=\rho(t)$. We then have $\frac{d f}{d t}=p X + q Y=(p + q \rho)X$.
Further $\int \frac{d x}{d t} \frac{\delta f}{\delta x} d t=-\int X Y d t$, that is given the right hand side =0, 
we have that $X \bot Y$. Further, when $X \bot Y$ implies $X^{*} \bot Y^{*}$, we have that 
$\xi \bot \eta$ gives an orthogonal base. Given $f^{\diamondsuit}$ harmonic or f linear, we have that $\frac{\delta^{2} f}{\delta x \delta y}=\frac{\delta^{2} f}{\delta y \delta x}$.

\subsection{Symmetric transformations}
Assume $f=\beta e^{\alpha}$, we then have $\frac{\delta}{\delta x} f/ \frac{\delta}{\delta y} f=
\frac{\frac{\delta}{\delta x} \log \beta + \frac{\delta \alpha}{\delta x}}{\frac{\delta}{\delta y} \log \beta + \frac{\delta \alpha}{\delta y}}$.
Given $\frac{\delta \alpha}{\delta x}=\frac{\delta \alpha}{\delta y}$, we have that over a contractible domain 
that $\frac{d y}{d x}=1$. Where $\alpha$ is symmetric, we have that  $\log \beta$ is symmetric.
When we are discussing a domain of symmetry $\{ \log \beta < \lambda \} \subset \subset \Omega$,
we regard projectivity $x \rightarrow y$ as necessary. Given $y=y(x)$ linear and $x \rightarrow y$ projective, 
We have that  $\log \beta$ is symmetric. 

U is radical if $U (\xi,\eta) \sim (\xi,\eta)$,
that is $d U \xi=d \xi$,$d U \eta = d \eta$. Note that given $X_{U}(\xi)=X_{U}(\eta)=0$, we have $X_{U}(f)=0$.
Assume $W f \sim \rho f$ in $L^{1}$. Given $d U=d U_{1} W=0$, we have that $\frac{\delta \rho}{\delta x}/\frac{\delta \rho}{\delta y}=const$, that is $\rho$
is symmetric in $x,y$. Assume $\{ \rho,f \}=X(f)$, then we have that
$<(\xi,\eta),(-Y,X)>=0$ iff $<\rho, (\frac{\delta^{2} F}{\delta x^{2}},\frac{\delta^{2} F}{\delta y^{2}})>=
\int \rho \Delta F d x d y=\int X_{2}(F)=0$. Thus, given $X_{2}(F)=0$, we have that $WF$ is harmonic. 
When $\{ \rho < \lambda \} \subset \subset \Omega$, then W defines a closed movement  
where it is defined. When $d U^{\bot}=\alpha d U$, where $\alpha$ is written i $u_{1},u_{2}$, the sub level sets
to $\alpha$ over $\Gamma$ obviously have cluster sets, that is are unbounded.

Assume $\frac{\alpha \xi}{\alpha \eta} \sim \frac{0}{0}$.
Let $\rho=\xi / \eta$, we then have $\frac{\frac{d}{d t} (\alpha \xi)}{\frac{d}{d t} (\alpha \eta)} \sim 
\frac{\frac{\alpha'}{\alpha} \rho + \frac{\xi'}{\xi} \rho}{\frac{\alpha'}{\alpha} + \frac{\eta'}{\eta}}$.
Assume $\xi \rightarrow {}^{t} \xi$ preserves character, we then have $<X(f),\varphi>=<f, \xi_{x} \varphi + \eta_{y} \varphi> + 
<f, X(\varphi)>$. Sufficient for this is that $f \bot (\xi_{x} + \eta_{y}) \varphi$. Assume $\frac{\delta G}{\delta y}=-\xi$,
$\frac{\delta G}{\delta x}=\eta$, then the condition is $-\frac{\delta^{2} G}{\delta x \delta y} + \frac{\delta^{2} G}{\delta y \delta x}=0$,
which is the case when G is linear in x,y or symmetric.  

Given $\{ \phi,\psi \}=0$, we have that symmetry for $\phi$ implies symmetry for $\psi$. Assume $\phi=U^{\bot} \psi$
and $\phi$ does not change sign on connected  components, for instance the boundary does not contain 
connected  components, we then have $\{ \phi^{2},\psi \}=0$ implies $\{ \phi,\psi \}=0$. In particular if $\phi^{2}$ 
is symmetric then $\psi$ is symmetric. 

\newtheorem{LemmaS}[def1]{Lemma (Symmetry is radical)}
\begin{LemmaS}
 Assume $\phi \in \mathcal{D}_{L^{1}}$ does not change sign on connected components, then relative
 involution, symmetry for $\phi$ is a radical property.
\end{LemmaS}

Assume $d U=\alpha d U_{1}$. Assume $\frac{d^{2}}{d t^{2}} U=\frac{d \alpha}{d t} \frac{d^{2} U_{1}}{d t^{2}} + \alpha \frac{d U_{1}}{d t} > 0$,
that is $\frac{d}{d t} \log \alpha > - \frac{d}{d t} \log \frac{d}{d t} U_{1}$. Given $(\xi_{1},\eta_{1}) \sim const$,
we have that $\frac{d}{d t} \alpha > 0$, that is $ \log \alpha$ is monotonous. Note that given $z \in \mbox{ supp }E_{1}=\Omega$,
we then have if we assume $\Omega=\log S$, where $S$ is formally a ball, then $e^{z} \in S$. 
U considered over S is algebraic. Note that if $X_{j}=d U_{j}$ is a reduced measure and $\Phi_{j}=X(f_{j})$, 
we have that $\Phi_{1}=\Phi_{2}$ implies $f_{1}=f_{2}$
(relative the topology for f). 

\subsection{Symmetry relative parameters}
$U \in$ center $\mathcal{G}$ implies $UV=VU$,  $\forall V \in \mathcal{G}$, this means that if $I \in \mathcal{G}$ 
we have that $UI=IU$.
Assume $d U=\rho d U_{1}$,$d U^{\bot}=\vartheta d U_{2}$ and $F(d U_{1},d U_{2})=F(d U_{2},d U_{1})$.
We then have $F(d U,d U^{\bot})=F(\rho d U_{1}, \vartheta d U_{2})$. Assume $=F(\vartheta \rho d U_{1},d U_{2})$
$=F(\rho \vartheta d U_{2},d U_{1})$. Then F is symmetric with respect to $d U,d U^{\bot}$.

Assume $F(u_{1},u_{2})$ symmetric through $\frac{\delta^{2} F}{\delta u_{1} \delta u_{2}}=\frac{\delta^{2} F}{\delta u_{2} \delta u_{1}}$.
Since $\frac{\delta F}{\delta u_{1}} \equiv 0$ implies the symmetry condition,
we see that one parameter movements satisfy the symmetry condition. Assume $\frac{\delta F}{\delta u_{1}}=-\tilde{Y}$,$ \frac{\delta F}{\delta u_{2}}=\tilde{X}$.
The regularity conditions can be given in the parameters to $u_{j}$, $j=1,2$. Assume $\tilde{F} \sim_{m} F$ (monotropy \cite{Dahn13}), where 
$\tilde{F}$ has compact translation and $\int_{\Gamma} d \tilde{F}=\int_{\tilde{\Gamma}} d F=0$,
where $\tilde{\Gamma}$ is the closed curve that corresponds to the consequent of $\Gamma$. 
We then have that $F$ is almost periodic. Assume $d U = \alpha d U_{1}$ and $d U^{\bot}=\beta d U_{2}$. 
When $<\alpha d U_{1}, \beta d U_{2}>=0$, given $\alpha,\beta$ algebraic (removable) 
we have that $<d U, d U^{\bot}>=0$. Assume $\alpha=e^{\phi}$ and $X(\phi)=0$, that is U analytic in the phase. We then have $X(\alpha f)=X(\alpha) f + \alpha X(f)$
$=X(\phi)(\alpha f) + \alpha X(f)$, that is $X(\alpha f)=\alpha X(f)$.

\subsection{Projectivity}
Assume $\mathcal{G} \times \mathcal{G}^{\bot}$ very regular. Assume $d U=\alpha d U_{1}$. Given $(I) \subset (I_{1})$, 
we have that $1/\alpha \rightarrow 0$ i $\infty$. In this case $U$ is not projective, when $U_{1}$ projective.
Assume $(I) \bigoplus (I^{\bot})=(I_{1}) \bigoplus (I_{2})$, Given $(I) \subset (I_{1})$ we must
have $(I_{2}) \subset (I^{\bot})$. Assume $\frac{d U^{\bot}}{d U_{2}}=\beta^{\bot}$. Given
$d U_{1} + d U_{2}=(1 + \gamma) d U_{1}=d x$, we have that $\gamma \rightarrow 0$ in $\infty$. 
Given $d U + d U^{\bot} = \alpha d U_{1} + \beta^{\bot} d U_{2}=(\alpha + \beta^{\bot} \gamma) d U_{1}=d x$,
we have that $\beta^{\bot} \sim (1 -\alpha) / \gamma$ in $\infty$. When $\beta^{\bot}$ is further bounded, 
we have that $\alpha \rightarrow 1$ in $\infty$! Given $\{ U,U^{\bot} \}=0$ we have $\xi \eta^{\bot} - \eta \xi^{\bot}=0$. 
Assume  $\xi^{\bot}=1/\xi$,$\eta^{\bot}=1/\eta$.
Thus, we have that $(\xi +1/\xi ) / (\eta + 1/\eta) \sim 0/0$, when $\xi,\eta \rightarrow 0,\infty$. 
When $ U + U^{\bot}$ is harmonic, $(\xi + \xi^{\bot})_{y}/(\eta + \eta^{\bot})_{x}=1$,
we thus have $(\xi + \xi^{\bot})/(\eta + \eta^{\bot}) \rightarrow 1$.
Assume in particular $\frac{\delta U^{\bot}}{\delta y}=\frac{\delta U}{\delta x}$,$\frac{\delta U^{\bot}}{\delta x}=-\frac{\delta U}{\delta y}$
that is $\xi^{2} + \eta^{2}=0$. We then have $\Delta U= \frac{\delta^{2} U^{\bot}}{\delta x \delta y} - \frac{\delta^{2} U^{\bot}}{\delta y \delta x}$.
\emph{ Thus, U is harmonic iff $U^{\bot}$ is symmetric. } 

Consider $f(\zeta_{T})=F(\gamma_{T})(\zeta)$. Given that U can be separated over f, we have not 
necessarily simultaneously separability over $\zeta$, $(f (\zeta_{T}) + f (\zeta_{S}))'=f' (\frac{d \zeta_{T}}{d T} + \frac{d \zeta_{S}}{d S})$.
Given $U_{T} + U_{S}$ projective, we have that not necessarily $f(\zeta_{T}) + f(\zeta_{S}) \rightarrow \{ \zeta_{T} \} \cup \{ \zeta_{S} \}$
continuous. On a contractible domain , where $w,w^{\diamondsuit}$ are harmonic, we have that $w = d f + d \overline{g}$
with $f,g$ analytic. Given $d w=d U + d U^{\bot}$, where the terms are harmonic, when d w is exact, 
then the movements can be chosen analytic.
\emph{Given $(d U,d U^{\bot}) \simeq (d U,0)$, that is the domain is a
translation domain, then the distributional contribution can be seen as negligible}.
Assume $d U=\alpha d U_{1}$. When $d U$ is locally reduced, we can  choose $d \alpha \neq 0$ outside 
the boundary.

\section{Conjugation using the Fourier transform}

\subsection{Orthogonal in the mean}
Note (\cite{Schwartz66}) that $(\dot{B})'$ is nuclear but not $(\mathcal{D}_{L^{1}})'$. 
When $\phi \in L^{1}$ implies $W \phi \in L^{1}$, we have that
$W^{\bot} \widehat{\phi} \rightarrow 0$, with $\widehat{\phi} \in \dot{B}$, that is $W^{\bot}$ 
has nuclear representation. When $W^{\bot}$ is projective, then $W^{\bot}=I-W$ is nuclear. 

Consider a subset of $\phi$ such that $W^{\bot}$ projective over $\widehat{\phi}$, say A. When 
we have that $W^{\bot}=I-W$ over A, then W is nuclear over A. Assume $A \bigoplus B=L^{1}$. 
Assume $(d W^{\bot})^{\bot}$ is given by $d V$. Assume $d V$ is not closed, but $\mid d V \mid^{2}$ 
closed. Let $\tilde{B}$ be a domain for $\mid d V(\varphi) \mid^{2}$. Given the domain for $ V^{2}$ 
is $B$ and $\mid d V^{2}(\varphi) \mid \leq \mid d V(\varphi) \mid^{2}$, we assume the difference is 
a zero-function.

\newtheorem{lemma8}[def1]{Lemma (Duality in the mean)}
\begin{lemma8}
 Assume existence of $V \in \mathcal{G}$ such that $M(V f)=\widehat{\frac{d f}{d x}}$ in $\mathcal{D}_{L^{1}}$.
 Then $U \rightarrow V$ satisfies $U^{\bot} M = M V U V^{-1}$.
\end{lemma8}

Solutions to partially differential equations constitute a closed subset of $\mathcal{D}'$, 
that is $\Phi_{1} E_{1} = \delta$ permits
annihilators to $\Phi_{1}$, that is $(d U)^{\bot} \sim d V$.
Sufficient to determine the character of the movement, is to consider action in the phase. Assume for this reason
$(I)=\{f \quad \log f \in L^{1} \}$. Assume $\log f = \phi$ and $W \phi \in L^{1}$. Starting from
$<\widehat{I}(W \phi) , \varphi>=0$ iff $<W \phi, \widehat{\varphi}>=0$. Assume $L^{1} = R(W) \bigoplus A$.
Define $(I)^{\bot}=\{ \varphi \quad \widehat{\varphi} \in A \}$. Given $W$ projective, $W^{\bot}=(I-W)$. 
Thus, we have existence of $V$ on $R(W)^{\bot}$, such that $\forall \phi \in D(W)$, there is a 
of $\varphi \in (I)^{\bot}$, such that 
$(I-W) \phi = V \widehat{\varphi}$. In the same manner  $<W \phi, V \widehat{\varphi}>=  0 = <W^{\bot} \widehat{\phi}, V^{\bot} \varphi>$
and given that $V$ or $W$ projective, we have that $W^{\bot} \simeq V$. Assume existence of G, such that
$W \phi =\{ G , \phi \}$, we then have $\{ F(G), \phi \}=\frac{d F}{d G} \{ G,\phi \}$, where 
we assume $\frac{d F}{d G}$ regular.

\section{Projectivity in graph norm}
 
\subsection{Completion}
Consider $Uf \in \mathcal{D}_{L^{1}}$ and $U^{\bot} \widehat{f} \in \dot{B}$ completed to $L^{1}$.
We then have $\mid U f - U^{\bot} \widehat{f} \mid^{2} \sim \mid U f \mid^{2} + \mid U^{\bot} \widehat{f} \mid^{2} - 
2 <U f, U^{\bot} \widehat{f}>$. Further, $\mid <U f, U^{\bot} \widehat{f}> \mid \leq \parallel U f \parallel \parallel U^{\bot} \widehat{f} \parallel$. 
Assume $d U^{\bot} \in L^{1}(d U)$, that is $dU^{\bot}/d U$ has isolated singularities. 
Consider $W=(U,U^{\bot})$, given $U^{\bot \bot} \simeq U$, we have that $W^{\bot} \simeq W$.
When $W$ is considered as a two parameter movement , we have that $d W$ is reduced implies $d U^{\bot} \prec \prec d U$. 
Note that if U is one parameter, we have that $U^{\bot}$ is not necessarily one parameter. 
\emph{When $d U^{\bot} = \alpha d U$
with $\alpha \in L^{1}(d U)$ this means that $\log \alpha \in L^{1}$ iff $\log \frac{1}{\alpha} \in L^{1}$}

Note that relative $<f,\widehat{g}>=0$, given $U f \bot U^{\bot} \widehat{g}$
implies $U^{\bot} \widehat{g}=0$, then projectivity for $U^{\bot}$ is
sufficient for symplecticity. 
Note that when $U^{\bot}=I-U$, we have that
$U^{\bot}$ has isolated zeros where $U$ have algebraic zeros.

We assume $T \Sigma = \{ \Delta U =0 \}$, that is the set where $U$ is harmonic. Assume $\widehat{\widehat{x^{j} f}} \rightarrow 0$ implies
$\mid x^{j} f + \rho \mid \leq c$, as $x \rightarrow \infty$. Modulo $C^{\infty}$ we can  assume f of real type,
that is $\mid \rho \mid \leq \epsilon$, Thus, $x^{j} f \sim 0$ in a neighborhood of $x_{0}$ and 
$\mid f \mid \leq c / \mid x \mid^{j}$, that is modulo zero sets, we have that $f \in B_{pp}$ if 
$\widehat{f} \in B_{pp}$. We assume $U \in \mathcal{D}_{L^{1}}'$
iff $U^{\bot} \in \mathcal{D}_{L^{1}}'$, given $T \Sigma = (T \Sigma) \cap (T \Sigma^{\bot}) \mid_{\mathcal{L}}$ (\cite{Dahn19})

Given $X_{j} {}^{t} X_{j} \sim X_{V} \sim {}^{t} X_{j} X_{j}$, we have  defined a normal operator.
When $X_{V} \sim (F,M)$ (\cite{Dahn19}), we have  further a normal model over the set where $\frac{x}{y} \rightarrow \frac{y}{x}$
projective. Over this set we have that $V=V^{\bot}$, that is the symmetry set defines a spiral domain.
Over $\Gamma$, we have that obviously $X_{V}$ defines a normal operator. However, we have that  
$d V$ is not BV, that is we do not have a determined tangent. Thus, given $X_{j}$ according to the above, 
the condition that a non-trivial M is not symmetric is necessary for a normal model with determined tangent.

\subsection{A normal model}
Assume $\log f \in \mathcal{D}_{L^{1}}$ with $\frac{d}{d x} \log f \in \mathcal{D}_{L^{1}}$.
Given $f = \vartheta \frac{d f}{d x}$, we have that when $\frac{1}{\vartheta} \in \mathcal{D}_{L^{1}}$,
$1/\vartheta \rightarrow 0$ in $\infty$. Given $f$ locally  polynomial,
we have that $\frac{d f}{d x} \prec f$, that is $\{ \log \frac{d f}{d x} < \lambda \} \supset \{ \log f < \lambda \}$
$\simeq \{ \log \frac{ d f}{d x} + \log \vartheta < \lambda \}$, that is $\vartheta$ can be chosen as mollifier. 
For instance $\Omega_{\epsilon}=\{ \log \frac{1}{\vartheta} > \epsilon \}$ and $\log f \in \mathcal{D}_{L^{1}}$ can 
be approximated with $H$, where $\Omega_{\epsilon} \rightarrow \Omega_{0}$. 

Assume $d V, d V^{\bot}$  BV measures such that $d V \bot d V^{\bot}$ implies $\Gamma=\{ 0 \}$. 
This is regarded as a normal model. 
 When V is a normal operator, the model is independent of 
 orientation. Given U a normal operator, we have that $\Gamma'=\{ U^{\bot \bot} = U \} \subset \Gamma$. 
 Given $(U U^{\bot})^{\bot} \bot U U^{\bot}$, we must have $\Gamma'=\{ 0 \}$, that is $U \bot U^{\bot \bot}$.

\subsection{Projectivity}
Given $f \rightarrow \parallel f \parallel_{G}$,
locally  1-1, then $(U,U^{\bot})$ can be considered as projective on the range and the dimension
for $f \in (I)$ is preserved. 
Assume $(U f)^{\bot} \simeq \{ V f \}$ and $U^{\bot} f \subset \{ V f \}$. Given $(U , U^{-1})$ is
projective, we can  write $V=(I-U)$.
Given $U + U^{-1} \equiv I$, we have that $(U,U^{2})$ is projective. 

  Sufficient for a desingularization is that
  $I \prec \prec dU$, that is invariant sets are zero sets. Note that $d U$ preserves compact sub level sets
  locally . Given $U,U^{\bot}$ ac, then $d U= d U^{\bot}=0$ implies that $U^{\bot}=U^{\bot \bot}$ 
  (U is then not projective) that is diagonal. 
  This means that $f \equiv 0$ (polar). 
  Assume $d U^{\bot}=\alpha d U$, where $U$ is defined through $U^{\bot}$.
  The condition $\frac{d}{d t} (U,U^{\bot}) f \neq 0$ does not exclude $\Gamma$. As long as U monotonous,
  a sufficient condition for $\Gamma=\emptyset$, is that $sgn \alpha$ is negative.
  Assume U projective, in the sense that $U^{\bot}=(I-U)$, then $UI=IU$ implies $\Gamma=\{ 0 \}$.
 
\section{Representation using reduced measures}

\subsection{Dependence of parameter}
Assume U a one parameter movement and $U f= \int f X(t) d t$, where t is a real parameter. Assume 
$\frac{d y}{d x}=\rho(t)$. Thus $X(f)=(\xi + \eta / \rho) \frac{\delta f }{\delta x}$. If $\vartheta=\frac{d x}{d t}$, 
we have that \emph{$X(t)=(\xi + \eta / \rho) \frac{1}{\vartheta} \frac{d f}{d t}$, which is a differential operator 
in a real parameter}. Thus, $X(t) d t$ gives a reduced measure. Assume $X,X^{\bot}$ define
one parameter movements and $X^{\bot}=\alpha X$. Linear independence means that, given X is 
projective, that $X=X^{\bot}$ implies $s=t=0$.
Otherwise, we have existence of $t_{0}$, such that $X(t_{0})=X^{\bot}(t_{0})$ 
through the boundary conditions.

\newtheorem{lemma13}[def1]{Lemma (Reduced measures)}
\begin{lemma13}
 Any one-parameter sequential movement U, where d U is BV, can be represented as a locally reduced measure outside the polar.
\end{lemma13}

Note that given $W$ a two parameter movement, such that $d W$ BV and $W^{N}$ is reduced, then $W^{N}$ can be 
written as one parameter. Assume $\frac{\delta U}{\delta x}=\xi$,
$\frac{\delta}{\delta x} U^{2} = 2  U \frac{\delta U}{\delta x}=2 U \xi \sim 2 \xi'$, 
that is we assume that U preserves
bi characteristics. Given $(\xi,\eta) \bot (-Y,X)$, we have that $U^{2} \simeq U$.

\subsection{Relative almost periodicity}
Assume $U * \alpha(f)=U(\beta * f)$, with $\alpha,\beta \in \mathcal{D}$. Given $f \in  B_{pp}$, we have that $M(f * \beta)=M(f)$. 
When $\alpha \rightarrow \delta$,
given $\int f* \beta d U=\int f d U$, can be seen as f pp relative dU.
Given $U I=I U$, we have that when $\beta \rightarrow I$, $Uf * \alpha \rightarrow UI f$. Note that 
I can be represented through translation over $\mathcal{D}_{L^{1}}'$.

 Consider $e^{i \lambda x} \widehat{f} \sim U^{\bot} \widehat{f}$ and $Uf \sim f_{t}$, 
 \emph{then almost periodicity means that U is normal somewhere}. Assume $U \in \mathcal{D}_{L^{1}}'$ close to the boundary
 with $\{ U f \} \subset \subset \Omega=nbhd \Gamma$. Through the condition on very regular boundary, there is 
 $W^{\bot} \in \mathcal{G}$, $W \neq U$, analytic close to the boundary. Assume further 
 $R(W^{\bot}) \subset C^{\infty}$. Thus, if $d UW=0$, we have when the movement is considered modulo $C^{\infty}$,
 that $A=I - W^{\bot} U^{\bot}=0$. Since $W^{\bot} U^{\bot}$ is compact, we have
 $N(A)=\{ 0 \}$, that is $d UW$ is projective on its range close to the boundary.
 Given $U f \in L^{1}$, as $U \rightarrow I$,
 we have that $U^{*} \widehat{f} \rightarrow 0$ in $\infty$, when $U^{*}$ is completed to 
 $U^{\bot}$ with preimage in $L^{1}$, we have that $U^{\bot} \widehat{f} \rightarrow 0$ i $\infty$, 
 but since $\dot{B}$ is not reflexive (cf. $\mathcal{D}_{L^{1}}$), it does not follow that 
 $U^{\bot \bot} \widehat{\widehat{f}} \rightarrow 0$ i $\infty$.

 Given f pp and $\widehat{f} \in B$, we have that $\widehat{f}$ pp. 
 In this case, there is a set G
 where $\sup \mid U_{1} \widehat{f} \mid$ is reached, that is we have a maximum principle through G.
 Close to $\Gamma$, G can be associated to general U. When $W$ is analytic, we have a regular approximation
 of $G$. Assume $W=V^{\bot}$, where $V V^{\bot}=V^{\bot} V$ with $V \neq V^{\bot}$.
 Given W is reflexive, we have that V is normal implies W is normal. We can in this manner 
 construct a normal system through $(W,W^{\bot}) \rightarrow (W,V) \rightarrow (V^{\bot} ,V)$.

 \subsection{Algebraicity}
Assume $Uf = e^{V \phi}$ a one parameter movement and $U_{S}$ a spiral approximation of $U$, 
such that $U_{S} f=\beta e^{V \phi}$. Consider U,V harmonic. Then $\log M$ is convex, where M 
is the maximum of U over an interval. Algebraicity for U means $U \simeq V$.

Note $\int d \widehat{Uf} = \int \xi \frac{d \widehat{f}}{d x} + \eta \frac{d \widehat{f}}{d y}$.
Thus, the closed property, exactness, the reduced property for the measure that defines U, is dependent of
the regularity for $\xi,\eta$. Given $\xi,\eta$ algebraic and $\Omega$ of positive dimension, then we have that $\int_{\Omega} d \widehat{ U(f)}=0$ 
implies $\int_{\Omega} d \widehat{f}=0$, that is the movement preserves pseudo convexity. The measure that 
is defined by $\xi d x - \eta d y$
can be considered as reduced, if $\frac{1}{\xi},\frac{1}{\eta} \rightarrow 0$ in $\infty$, 
that is $I \prec \prec U$. 

Assume $d U^{N}$ reduced, which implies $d U$ BV. Assume $d U^{\bot} \simeq g^{*} d x^{*}$,
where $g^{*}$ polynomial . Given $g$ is reduced, there is $g \rightarrow  1/g$
projective. As long as the completion of $g^{*}$ is algebraic, we can choose $d U^{\bot}$ as BV.

\section{Spiral approximations}
 Assume $U_{S} \rightarrow U$ where $U_{S}$
has support in the polar and that $\int \mid d U \mid < \infty$, that is the complement to the polar 
has finite D-integral, then the spiral closure can be determined. 
Determination of spiral closures is a partially unsolved problem, given the complement to the range 
has infinite D-integral. Given $<U_{S} f, U_{S} \widehat{g}> \simeq <U_{S}^{2} f,\widehat{g}>$,
when $U_{S}^{2} \rightarrow U$, then $U$ has two possible limits. 

Assume $(\sqrt{U_{1}} + i \sqrt{U_{2}})^{2}=U_{1} - U_{2} + 2 i \sqrt{U_{1} U_{2}}$. Assume for instance
$U_{1} \varphi \in L^{1} \cap L^{2}$.
Consider the change of variables $ \psi : (x,y) \rightarrow (u,v)$ (sequential movements). When u 
is fixed, the spiral approximates a circle in a neighborhood of u=const. When v is fixed, 
the spiral approximates a line in a neighborhood of v=const. Assume $\tilde{F} (u,v)= F \circ \psi(x,y)$ 
and $d \tilde{F}=\alpha d F$, 
where $\alpha$ is regular outside $u=v$, that is
given $d F(x,y)$ regular, we assume $d \tilde{F}(u,v)$ very regular. ' 

\subsection{A separation property}
Concerning iteration $X^{2}(f)=\xi^{2} \frac{\delta^{2} f}{\delta x^{2}} + \eta^{2} \frac{\delta^{2} f}{\delta y^{2}} + X(\xi)
\frac{\delta f}{\delta x} + X(\eta) \frac{\delta f}{\delta y} + \xi \eta (\frac{\delta^{2} f}{\delta x \delta y} + \frac{\delta^{2} f}{\delta y \delta x})$. Given $X$ corresponds to $U$, such that 
$X(\xi)=X(\eta)=0$ and $\frac{\delta^{2} f}{\delta x \delta y}=-\frac{\delta^{2} f}{\delta y \delta x}$, we have that $X^{2}=\xi^{2} \frac{\delta^{2} f}{\delta x^{2}} + \eta^{2} \frac{\delta^{2} f}{\delta y^{2}}$.
Given $X^{2}=0$, we have that in this case $\frac{\delta^{2} f}{\delta x^{2}} / \frac{\delta^{2} f}{\delta y^{2}}= - \eta^{2} / \xi^{2}$.
Note that for instance $\xi^{2}$ polynomial or real
does not imply the same property for $\xi$. 
Assume $\frac{\delta g}{\delta x}= \frac{\delta^{2} f}{\delta x^{2}}$
and in the same manner for $\delta / \delta y$. Let $X''(g)=\xi^{2} \frac{\delta g}{\delta x} + \eta^{2} \frac{\delta g}{\delta y}$
and $X'(f)=X(\xi) \frac{\delta f}{\delta x} + X(\eta) \frac{\delta f}{\delta y}$.
Choose $h$ such that 
$\frac{\delta h}{\delta x}=\frac{\delta^{2} f}{\delta x \delta y}$ and $\frac{\delta h}{\delta y}=-\frac{\delta^{2} f}{\delta y \delta x}$
For this reason choose h symmetric, that is such that $\frac{\delta h}{\delta x}=\frac{\delta h}{\delta y}$. 
\emph{We can then solve $X^{2}(f)=Y(h)$ through $X''(g) + X'(f)=Y(h)$}. 

Assume $\frac{\delta Z}{\delta x}=X(\xi)$ and $\frac{\delta Z}{\delta y}=X(\eta)$.
We then have, when $X(f)=0$, that the condition $X' = 0$, is equivalent with $\xi X(\eta) - \eta X(\xi) = 0$.
For $X''=0$, we have to assume $\eta^{2}/\xi^{2} \simeq \eta / \xi$ over $g$.

\newtheorem{prop13}[def1]{Proposition (A separation property)}
\begin{prop13}
 Assume existence of $h$ such that $X^{2}(f)=Y(h)$, such that Y is exact over h. Then there
 is $g$ such that $Y(h)=X'(f) + X''(g)$. When $X(\eta) / X(\xi) \simeq \eta^{2} / \xi^{2} \simeq \eta / \xi$,
 the movements to the right hand terms can be chosen analytic.
\end{prop13}

Consider the extended domain in $(x,y) \rightarrow (x,\frac{y}{x};y, \frac{x}{y})$.
Assume $\frac{y}{x} \rightarrow \frac{x}{y}$ projective with $\frac{d y}{d x} \neq 0 \Rightarrow \frac{d x}{d y} \neq 0$
that is $y(x) \rightarrow x(y)$ regular. Given f holomorphic, we have that $\min f= \max 1/f$. In the extended plane
we assume minimum (maximal domain ) are symmetry points, that is $\Omega \ni (\frac{x}{y},\frac{y}{x})$ iff
$(\frac{y}{x},\frac{x}{y}) \in \Omega$. Given a maximum principle, there is no point outside the symmetry
set,  where minimum (spiral) is reached. Note that symmetry implies $(u,v) \rightarrow (v,u)$ projective on the domain. 
Note that if we consider f modulo $C^{\infty}$ with  $f={}^{t} f$, then the symmetry property is radical (\cite{Dahn15}), 

\subsection{The parameter space}
Consider $(Ug,U^{\bot} g)=g(U,U^{\bot})$. Assume for the iterated symbol, $g_{N}(U,U^{\bot})=0$ implies $s=t=0$, but we can have
$g(U,U^{\bot})=0$, as $s,t \rightarrow \infty$.

\newtheorem{lemma16}[def1]{Lemma (The spiral in cylindrical parameters)}
\begin{lemma16}
 Assume $(x,y) \rightarrow (u_{1},u_{2})$ regular, then the equation for the spiral can
 be written as $d U_{1}(f)(u_{1},u_{2})$.
\end{lemma16}

We consider $U_{1} \gamma=\gamma(u_{1},0)$ and $U_{2} \gamma = \gamma(0,u_{2})$. Further, 
given $d y / d x=\rho$ and $\xi',\eta'$ are coefficients to Y(f) according to the above, 
we then have that the coefficients to $\xi,\eta$ in $X(f)$ are given by 
$\xi=\xi' \frac{\delta x}{\delta u_{1}} + \eta' \frac{\delta x}{\delta u_{2}}$
and $\eta=\rho \xi'$.

Assume $Y(f)=\xi'(u_{1},u_{2}) \frac{\delta f}{\delta u_{1}} + \eta'(u_{1},u_{2}) \frac{\delta f}{\delta u_{2}}$.
Given the movements in $\mathcal{G}$ are taken in sequence,
we have that  $Y(f)=0$. In particular given $f \equiv const$ over a spiral axes, we have that 
$Y(f)=0$ over the spiral. The condition $f=\widehat{f}$ on $\Sigma$ depends of the topology.
Note $\log(regel) \simeq translation$ (\cite{Lie96})
is not necessarily unique, that is $d U^{\bot} = \rho d U_{1}$ does not uniquely determine 
$U^{\bot}$ outside the boundary.

\subsection{The spiral as a measure zero set}
Note that when $U^{-1}$ is projective, we have that $U^{-1} d V=0$ implies $d V=0$, that is closed 
forms are mapped on closed forms, which is necessary for a normal model. When $U f(x)=f ({}^{t} U x)$ 
and $U \simeq {}^{t} U$ the condition is that $d {}^{t} U x = d x$. For instance when $U \simeq {}^{t} U$ 
and $U$ is translation, we have that $f(x) + c \simeq f(x + \eta)$, with $\frac{d \eta}{d x}=0$
and given $\eta$ ac, we have that $\eta = const$.
Note that if $U U^{\bot}=0$ and $U^{2} \simeq U$ (projectivity), we have that $(U + U^{\bot})(U-U^{\bot}) \simeq
U - U^{\bot}=I$. Note that if $U^{\bot}=(U-I)$, we have that $U^{\bot} U \simeq U - U^{2} \simeq 0$. 

Given convexity we have that
$f(1 - \frac{d y}{d x}) \leq \frac{1}{2}(f(x) + f(y))$. Given $\frac{d y}{d x}$ is projective, $f(\frac{d y}{d x}^{\bot}) \leq f(x) - f(y)$.
\emph{ In particular when $d U^{\bot} \bot d U$, $f (d U^{\bot}) \leq f(I) - f(U)$ } and for instance 
$f((1 - \alpha) d U) \leq f(U) - f (U^{\bot})$ with $d U \neq 0$ and $d U^{\bot}=\alpha d U$.

Assume that we have existence of G, such that $\int \frac{d G (\gamma)}{d \gamma}d \gamma \leq 
M(U \gamma) - M(U^{\bot} \gamma)$. Note $\{ G(\gamma),\gamma \}=\frac{d G}{d \gamma} \{ \gamma, \gamma \} \equiv 0 $
is independent of $\gamma$, Thus, $(I_{S}) \subset (I_{G})$, that is given G defines a measure and $\gamma$
a zero function to G, we have that the spiral defines zero functions. When the right hand terms 
are analytic, G is however not necessarily analytic.

\newtheorem{prop7}[def1]{Lemma (The spiral as zero sets to a measure)}
\begin{prop7}
 Assume existence of G such that $\int \frac{d G (\gamma)}{d \gamma}d \gamma \leq 
M(U \gamma) - M(U^{\bot} \gamma)$. Given that G defines a measure, the spiral is 
in the zero set corresponding to G. 
\end{prop7}

Given $N(I_{G})$ not removable,
we have that the same holds for $N(I_{S})$. Assume $\dim U \leq \dim U^{\bot} + \dim G$.
Given $G \gamma=0$ implies $\gamma \in (I_{S})$, we have that $(I_{S})$ is maximal in the sense
that there is not a continuation of $U_{S}$, that is ``maximal spiral rank''. Assume $(I-G)$ analytic on 
compact sets, then the dimension
outside the boundary, that is over regular points, is locally given by analytic functions. 
However given $(I_{S})$
is defined through $G \gamma=0$, where $U$ projective, then G is not projective. 
$(I_{S})$ can be seen as a subset of the complement to an analytic set. 

Assume $Y/X \sim_{0} 0 / 0$ implies that
$Y_{x} / X_{y} \neq const$ (vorticity), that is we assume that the movement does not change character 
and orientation simultaneously close to the boundary. Assume the movement is one parameter and 
$U_{1}$ is locally  reduced, that is $U_{1}f=0$ implies $f=0$,
we then have $\delta U / \delta t \neq 0$. The condition $d U$ is of order 0,
implies $\delta U / \delta U_{1},\delta U / \delta U_{2}$ constant close to the boundary $\Gamma$. 
We thus assume
that a neighborhood of the boundary can be generated by translation and rotation. 

Assume $\delta U / \delta U_{1}=\alpha,\delta U / \delta U_{2}=\beta$ constants close to the boundary. We then have 
$\Delta_{u_{1},u_{2}} U =0$, that is we can assume that close to the boundary, that $U$ is harmonic 
relative $u_{1},u_{2}$, why the boundary is of order 0.
Further, $\frac{\delta^{2}}{\delta u_{1} \delta u_{2}} U^{\bot}=\frac{\delta^{2}}{\delta u_{2} \delta u_{1}} U^{\bot}$
that is $U^{\bot}$ is ``symmetric'' close to the boundary. 
In the same manner , assume $d U = \alpha d U_{1}$
where $\alpha$ linear i 1,x,y. we then have $\frac{\delta^{2} \alpha}{\delta x \delta y} = \frac{\delta^{2} \alpha}{\delta y \delta x}$
that is $\alpha$ is symmetric i $x,y$.

\subsection{The mean spiral}

Consider $M(d U) \leq M(U) - M(I)$, Given $U$ projective and ac we have that the left hand side $\sim_{0}$ 
the right hand side. Further, $M(d U) - M(d U^{\bot}) \leq M(U) - M(U^{\bot})$, given $d U - d U^{\bot}=d I$.
Thus, $M(d I)= I \leq M(U-U^{\bot})$, that is given $d U$ projective (in the mean) there is not space for a spiral.

\newtheorem{lemma17}[def1]{Lemma (Projectivity in the mean)}
\begin{lemma17}
 When $d U$ is projective in then mean, it does not define a spiral.
\end{lemma17}

The max-principle for bi linear forms gives $\mid (U f, U^{\bot} \widehat{g}) \mid \leq C \parallel U f \parallel_{1}^{1/2}
\parallel U^{\bot} \widehat{g} \parallel_{1}^{1/2}$, for a constant C. Given $U$ is considered as 
a normal operator, we can define $\Gamma=\{ U=U^{\bot} \}$. 
Assume $\eta_{x}=\frac{\delta G}{\delta y}$ and $\xi_{y}=\frac{\delta G}{\delta x}$ 
and $\{ G,f \}=0$. That is the movement  corresponding to $\eta_{x},\xi_{y}$ constants, is
sequential, we could say that \emph{$U_{S}$ is sequential in the mean}. Note for $d U_{s}$,
we have $\frac{d \eta}{d \xi}=\frac{-d x + \kappa d y}{\kappa d x + d y}$ and $\eta$ can be completed
in $L^{1}$ to $d \xi \bot d \eta$. Note where $d U^{\bot} / dU \leq 1$, $\Gamma$ is $d U^{\bot} / d U=1$.

\subsection{Partial transforms}
Consider the problem when $U_{1}^{\bot} f= U_{2}^{\bot} f$ implies $U_{1}=U_{2}$. 
Consider $S(f)=x \frac{\delta f}{\delta x}$. Note that ord S(f)= ord f. Using the partial
Fourier transform, according to $\mathcal{F} (y \frac{\delta f}{\delta x}) \sim -x^{*} \frac{\delta \widehat{f}}{\delta y^{*}}$,
(\cite{Dahn19}) that is given $G(y,x)=y \frac{\delta f}{\delta x}$, we have that ${}^{t} \widehat{G} \sim \widehat{G}$. 
Given $S(f)=0$, we have that $\mathcal{F} S(f)=\widehat{f}$ iff $S(\widehat{f})=0$. Further, given $S(f) \in B$
we have that $\frac{\delta f}{\delta x} \in \dot{B}$. G(x,x)=0 implies $\widehat{G}(x^{*},x^{*})=0$ and conversely. Let 
$-T(f)=\mathcal{F}S(f)$, we then have $-T(f)=\widehat{f} + S(\widehat{f})$, $S(\widehat{f})=\widehat{S}(f) - \widehat{f}$
and given $S$ projective, we can write $S(\widehat{f})=\widehat{S}^{\bot}(f)$. Assume $\exists g \in B_{pp}$ such that $g = \frac{\delta f}{\delta x}$, that is $d f=g d x$ and
$d f=0$ implies $g=0$. Given T projective $B_{pp} \rightarrow \widehat{(I)}$, we can  write
$\widehat{S} \bigoplus \widehat{S}^{\bot}(f)=I \widehat{f}$. The condition $\widehat{S}^{\bot}(f)=0$ implies $\widehat{f}=0$
on compact sets, can be interpreted such that $I \prec S$. Consider $T_{\rho}(f)= \rho f$ where $\rho$ is complete and symmetric. Given $T_{\rho}$ projective and symmetric,
it gives an orthogonal base for $(I)$. In particular when the polar is defined through continuation  $T_{\rho}(\xi,\eta)=\rho (\xi,\eta)$
(symmetric) the polar can be chosen as orthogonal, that is $\xi' / \eta' \sim \eta / \xi$.
Assume for instance
$\xi = \frac{\delta G}{\delta y}$ and $\eta=-\frac{\delta G}{\delta x}$. we then have $X(f)=\{ G,f \}=0$,
gives that a set of symmetry for f has a corresponding set of symmetry for $G$. Assume $\rho(x,y)$ symmetric in $(x,y)$. 
Where $\rho \neq 0$, $T_{\rho}$ is locally 1-1.
Assume $\Sigma$ a domain for $f$ and $\Sigma_{\rho}$ a symmetric continuation. Given $U U^{\bot}=U^{\bot} U$ and $d U^{\bot}=\alpha d U$
we have that $\alpha \xi = \frac{1}{\alpha} \xi^{\bot}$ and $\alpha \eta = \frac{1}{\alpha} \eta^{\bot}$.
Given $\alpha$ symmetric, we have that  U is symmetric iff $U^{\bot}$ is symmetric, where 
$\alpha \neq 0$. The condition $d U \rightarrow d I$ regularly corresponds to a contractible domain . 
When we assume F symmetric
in a neighborhood of $d I$, a neighborhood of $d I$ can be generated by $d U$. Given $d U=0$ (analytic) we have that
for $(x,y)$ such that $d U \rightarrow d I$ regularly, that $(x,y) \in \mbox{ ker } F$. Consider $d U \rightarrow d U^{\bot}$
through $\rho$ symmetric. Then there is where $d U$ is projective, an orthogonal base. In particular 
$(d U)^{\bot} \simeq d U^{\bot}$.

\section{Maximal rank at the boundary}
The fundamental problem, we consider is a movement (the spiral) that is assumed to have maximal rank 
at the boundary (cylindrical). In presence of a maximum principle, we do not have hypoellipticity 
because of presence of trace. We assume the boundary very regular, that is we have existence of 
a regular approximation $V$ (not unique) of $\Gamma=\{ U=U^{\bot} \}$. Conversely we assume that the range for 
U always have points in common with some regular V. Given the movement is defined by $d V$ BV, 
we have determined tangent, that is V is not a spiral.

\subsection{Geometry}

In hyperbolic geometry, the orthogonal to a
hyperbolic movement  (euclidean translation) are parallel planes (\cite{Riesz43}). Any multi valent 
surface given by a holomorphic function,  
can be approximated through Puiseux expansions (\cite{Oka60}) $U_{s}^{\bot} \rightarrow U^{\bot}$. Starting with 
a multi valent surface, the cylinder web is regarded as a boundary.
Individual leaves intersect the web, not necessarily in a point. Assume $u_{1},u_{2}$ coordinates 
for the cylinder with 
$u_{2} / u_{1}$ constant, then the transform $(u_{1},u_{2}) \rightarrow (u_{1},u_{2}/u_{1})$ (ähnlich transform (\cite{Lie91})),
that is a (monotonous) one parameter movement . Consider $d U_{S}(x,y)=d U_{1}(u_{1},u_{2})$  $\sim d U_{1}'(u_{1},u_{2}/u_{1})$ 
as monotonous.

We assume $U$ factorized over
$U_{1},U_{2}$. When U is defined by $\Delta U =0$ we have that $U \sim {}^{t} U$. 

\newtheorem{prop5}[def1]{Lemma (Translation is a movement of real type)}
\begin{prop5}
 Assume $U$ factorized over $U_{1},U_{2}$ and $\mid U f \mid \sim e^{\mid w \mid \mid \phi \mid}$,
 with $w$ regular ($d w(\phi) \simeq \rho d \phi$ and $\rho$ regular).
 When $U$ has real type, the corresponding domain is a translation domain.
\end{prop5}
(cf. \cite{Lie96}, Ch. 9)
When $d U=d U_{S}$ is 1-dimensional in 3-space
and $d V \bot (d U, d U^{\bot})$, the maximal rank for V must be 1.

\subsection{Maximal rank on a subset}

Assume $d U^{\bot} = \alpha d U$, where $d U^{\bot}$, is seen as continuation of d U.
Given the continuation not analytic, the domain for d U analytic is maximal. 
\emph{When $U \simeq {}^{t} U$, we have that $U \rightarrow {}^{t} U $ is projective on the range.}

\newtheorem{prop6}[def1]{Lemma (The spiral is projective on its range)}
\begin{prop6}
 Assume the spiral is defined by $U \simeq {}^{t} U$ relative $<f,\widehat{g}>$, where $f,g \in \mathcal{D}_{L^{1}}$.
 Then the $U \rightarrow {}^{t} U$ is projective on its range.
\end{prop6}
We assume $U \rightarrow U^{\bot}$ is defined by $X(f) \rightarrow X^{\bot}(\widehat{f})$.
When the graph norm $\parallel f \parallel_{G}=0$, we have that $X(f)=0$ implies $X^{\bot}(\widehat{f})=0$. 
Further given $X \bot X^{\bot}$ and $\parallel f \parallel_{G}=0$, we have
$X=0$ iff $X^{\bot}=0$. Invariant sets corresponds in graph norm to $(U,U^{\bot})=(I,0)$. 
Assume $(d U)^{\bot} \simeq d V$. When $U$ is projective, we thus have that $d V \simeq d U^{\bot}$. 

Given $d U = d U^{\bot}$, then $d U \rightarrow d U^{\bot}$ is projective on $\mathcal{D}_{L^{1}}(\Gamma)$.
Every point on the cylinder web M, can be reached by spirals and given M, $u_{1}/u_{2}$ uniquely defines a spiral as
a one parameter curve.

\subsection{Weights}

 Assume $\mid U f \mid \sim \mid w \mid \mid f \mid$, where $w(tx,ty)=t^{1/2} w(x,y)$, $w(x,y) \neq 0$
 we then have $w'(t) \sim t^{-1/2}$, that is $1/w \rightarrow 0$ and $1/w' \rightarrow \infty$, as $t \rightarrow \infty$.
 Further $w'' \rightarrow 0$ and $w'' / w \rightarrow 0$, as $t \rightarrow \infty$. For a reduced measure 
 we have that $w^{(k)} / w \rightarrow 0$ in $\infty$. Note starting from $\mid W e^{\phi} \mid \leq e^{\mid w \mid \mid \phi \mid}$,
 where w does not change sign (on connected  components $\ni \infty$), given $w$ convex, we can  assume that 
 $w'' / w \rightarrow 0$ in $\infty$ one sided.
 
 Assume
 $U_{1},U_{2}$ is taken sequentially, with $\mid \widehat{U_{2} f} \mid \leq \mid w_{2} \mid \mid \widehat{f} \mid$
 and in the same manner  for $U_{1}$ with $w_{1}$. Given $w_{2} w_{1}$ complete with 
 $\parallel \mathcal{F}^{-1}(w_{2} w_{1} \widehat{f}) \parallel_{p} \leq \parallel f \parallel_{p}$,
 $w_{2} w_{1}$ can be considered as a $L^{p}-$ multiplier (\cite{Schechter70}).
 
Assume $\xi_{2}=2 f \xi$ and $\eta_{2}=2 f \eta$. We then have $(\eta_{2})_{x} - (\xi_{2})_{y}=2 f
(\eta_{x} - \xi_{y}) - 2 X^{\diamondsuit}(f)$. Thus, we have that $X(f^{2})$ is harmonic if $X(f)$ 
is harmonic and $X^{\diamondsuit}(f)$ is analytic. Further $X(f^{2})=2 f X(f)$. To determine dimension
for $Uf$ over $f \neq 0$, it is thus sufficient to consider $U f^{2}$.

 \subsection{Two sided limits}
Note that the two mirror model (surjective in the plane (\cite{Dahn14})) assumes two sided limits.
Assume $Vg \bot U f$ and $U^{-1} g=f$, for some $g \in R(U)$, which is implied by $U^{-1}$ surjective on
$R(U)$. Then when $V \not\equiv 0$, we can assume $V$ locally reduced. 
For instance $\lim_{\eta \rightarrow 0} (U_{1} f + U^{-1}_{1} f)=f(x + 0) + f(x - 0)$,
that is \emph{maximal rank can be compared with $\lim_{\eta \rightarrow 0} f(x + \eta)=\lim_{\eta \rightarrow 0} f(x - \eta)$},
that is $U + U^{-1}$ has maximal rank if the corresponding defect indexes are equal and zero.

\newtheorem{prop9}[def1]{Lemma (A two sided limit has maximal rank)}
\begin{prop9}
 Assume $d U = \alpha d U_{1}$, where $U_{1}$ has maximal rank and $\alpha \rightarrow 1$ regularly as $\eta \rightarrow 0$.
 Assume $d V=-\alpha d U_{1}$.
 Then $U + V$ has maximal rank, when $\lim_{\eta \rightarrow 0} U f=\lim_{\eta \rightarrow 0} V f$.
\end{prop9}

Assume $R^{2}=R_{1}R_{2}=R_{2}R_{1}$ projective.
Given $R^{2} \simeq R$, then R is projective, but not unique!
Given $R_{1} \rightarrow R_{2}$
projective, we have that for the corresponding Cayley indexes, $d_{-}=d_{+}$. Assume $(f,g)=<f, \widehat{g}>$ and $(U f, g)=(f , U g)$. 
Given U projective in graph norm with defect indexes 0,
we have maximal rank.

Assume $f=\beta e^{\phi}$ and $\beta=e^{\alpha}$, we then have $d U^{2} \sim U d U(f) \sim \beta (X(\alpha) + X(\phi))e^{\phi}$.
Assume $\beta=\beta(\phi)$, we then have $X(f) / f = X(\phi) (\frac{d \beta}{d \phi} + \beta)$, where 
we assume $\frac{d \beta}{d \phi} \neq 0$.  

\newtheorem{lemma1}[def1]{Lemma (Factorization lemma)}
\begin{lemma1}
 When $d U(f)=\{ G , f \}$, for some G, $U$ can be factorized into
 $d U (\beta e^{\phi})=d V (\beta) e^{\phi} +  d W(\phi) e^{\phi}$.
\end{lemma1}
that is $e^{-\phi} d U(f) = \{ G,\beta \} + \beta \{ G, \phi \}$

\subsection{Factorization}
Consider $Y(f)=\alpha \frac{\delta f}{\delta u_{1}} + \beta \frac{\delta f}{\delta u_{2}}$. Thus
if V is a movement defined by Y, we have $\alpha=\frac{\delta V}{\delta u_{1}}$ and $\beta=\frac{\delta V}{\delta u_{2}}$.
In particular when $V=U_{1}$, we have that $\beta=0$. Further for instance $\alpha=\xi \frac{\delta u_{1}}{\delta x}$ and 
$\beta=\eta \frac{\delta u_{2}}{\delta y}$. We then have $Y(f) \simeq X(f)$.

  Concerning $Y(f)=0$, given $U$ can be factorized through $V(u_{1}) W(u_{2})$ (sequential movement ),
  we then have $\beta_{u_{1}} - \alpha_{u_{2}}=0$. That is sequential movements are harmonic with respect 
  to the cylinder. Assume $d U = \rho d I$, where $\rho \rightarrow 1$
  regularly, when $U \rightarrow I$. When $\alpha,\beta$ are constants, we consider 
  the representation as not contractible relative $u_{1},u_{2}$. When $\{ \beta=const \}$ is locally negligible,
  the domain is a translation domain. Consider $\frac{\delta f}{\delta u_{1}}=\frac{\delta f}{\delta x} \frac{\delta x}{\delta u_{1}} + 
  \frac{\delta f}{\delta y} \frac{\delta y}{\delta u_{1}}$. Given $x,y$ linear in $u_{1}$, we have thus
  that $\frac{\delta f}{\delta u_{1}} \sim X_{U_{1}}(f)$. 

  \newtheorem{lemma2}[def1]{Lemma (Second factorization lemma)}
  \begin{lemma2}
   Consider $Y(f)=\alpha \frac{\delta f}{\delta u_{1}} + \beta \frac{\delta f}{\delta u_{2}}$.
   When $Y(f)=0$ defines a sequential movement in $u_{1},u_{2}$, we have that Y defines a harmonic movement
   (on the cylinder). 
  \end{lemma2}
 
Assume $U \frac{d}{d x} F=\frac{d}{d x} V F$ and when $F=e^{\phi}$, we have that $VF=e^{W \phi}$. 
Given $U$ algebraic,
we then have that $U \phi \sim W \phi$. Assume $U U^{\bot}=U^{\bot} U$ with U
projective $(U + U^{\bot})=I$. If further $I \frac{d}{d x} f = \frac{d}{d x} I f$, we have that 
$V V^{\bot} = V^{\bot} V$. Given $U^{2}=WU$, U is projective where $W=I$.
Given $U^{2}=U+V$, U is projective where V=0. 

Note that given $U^{\bot}=0$ implies $U =I \neq 0$, we have that
$N(U^{\bot}) \subset R(U)$.
For instance $U=WU_{1}$, where $d W \rightarrow d I$ with $\mbox{ ker} W=\{ 0 \}$, we then have that
$N(U_{1}) \subset R(W)^{\bot}$. If $WU_{1}=U_{1} W=0$, that is $W \bot U_{1}$, given $W$ projective, 
we have that $U_{1} \simeq W^{\bot}$ and $W W^{\bot} \simeq W - W^{2} \simeq 0$. When W is projective 
we have that $WU=0$ implies $U=0$ and UW=WU implies that U is locally 1-1. Note that when W is an 
extension and in the weak sense $UW=WU$, that is $W \rightarrow {}^{t} W$
preserves character, then $W$ must be projective if U is projective.

Consider as a max-principle, $F(z) \rightarrow F(P(z)) \rightarrow F (e^{w})$ (\cite{Nishino75-1})
that is $\max \mid F \mid < \infty$, where $z \in \Omega$, where $\Omega$ is algebraic or exponential.
When $F(z)=F(R z)$, where $Rz=e^{w}$ a restriction, we have that $F(z) \simeq \widehat{F}(w)$. 
Determine $\Sigma=\{ w=G(z) \}$ such that $F$ is reduced over
$\Sigma$, for instance $\mid G \mid \leq \mid F^{N} \mid$. For instance $\widehat{F}(w)=F(G(z)) \rightarrow \infty$, 
when $G(z) \rightarrow \infty$ as $z \rightarrow \infty$. 

Consider $\Omega$ a cylinder and assume $w=P(z) \in \Omega$ (algebraic polyhedron). A max-principle
is $d U$ BV over $e^{\Omega}$. Consider $d \widehat{U}(\Omega) = d U (e^{\Omega})$. The condition
$d \widehat{ I U}$ BV over $\Omega$ means that $d U$ is BV over $\Omega$ and $e^{ U \Omega}$ finite.

\subsection{Desingularization}
Note 
$\int_{\Gamma} \xi d x + \eta d y = \int_{(\Gamma)} (\eta_{x} - \xi_{y}) d x d y$. 
Given $\xi,\eta$ locally  bounded, we have over compact sets that we have a finite D-integral,
that is d U BV locally. Consider $\int_{\Omega} d U - A \int_{\Omega} d I \simeq \int_{\Omega} d U^{\bot}$.
In this case the projectivity for $U$ is dependent of the domain $\Omega$. For instance 
$\mid \frac{\delta f}{\delta x} \mid^{2} + \mid \frac{\delta f}{\delta y} \mid^{2} \geq c \frac{\mid x \mid^{2} + \mid y \mid^{2}}{\mid \xi \mid^{2}}$,
when $\mid x \mid,\mid y \mid \rightarrow \infty$ and $\mid \xi \mid \leq \mid \eta \mid$. Thus, given $0 \neq \mid \xi \mid$ locally  bounded
we have that $\{ d U \leq A \} \subset \subset \Omega$. Assume $\xi$ is independent of $\mid y \mid$, 
then we have that
$\mid y \mid / \mid \xi \mid \rightarrow \infty$, as $\mid y \mid \rightarrow \infty$. In the same manner 
when $\eta$ is independent of $\mid x \mid$. For instance $\eta_{x}=\xi_{y}=0$. This can be seen as a
``desingularization''. 

Given $dU=\alpha d U_{1}$ and $\alpha (x) \rightarrow x$ is
bounded, we have that the sub level sets are relatively compact, that is $\alpha$ reduced. Outside 
the polar to $U_{1}$, we have that $\alpha \equiv 0$ implies $\xi=\eta=0$. When $\xi,\eta$ are reduced, 
then $\alpha$ reduced.

\subsection{Projectivity}
Obviously, \emph{$U \in \mathcal{G}$ does not imply that U projective}.
When $U \in \mathcal{G}_{ac}$, $U' \in \mathcal{G}$, we have that $U U' \in \mathcal{G}$, 
that is $\frac{d}{d t} U^{2} \in \mathcal{G}$. 
Given $U,U' \in \mathcal{G}_{1}^{(1)}$
we have that $U^{2} \in \mathcal{G}_{1}^{(1)}$, but the converse assumes a regular inverse. Further 
$U \in \mathcal{G}_{1}$ and $\log U \in \mathcal{G}_{1}^{(1)}$ implies $ U U' U^{-1} \in \mathcal{G}_{1}$, 
that is $U' \in \mathcal{G}_{1}$.
Assume $\mathcal{G}_{1} \bigoplus \mathcal{G}_{1}^{\bot}$ has maximal rank. We have that
$U^{N} \in \mathcal{G}_{reg}$ (can be approximated regularly, that is contractible) does not imply 
$U \in \mathcal{G}_{reg}$. However, given $U$ reduced (has no invariant sets), we have that 
$U^{2} \in \mathcal{G}_{reg}$ implies $U \in \mathcal{G}_{reg}$.
Thus, when $\mathcal{G}$ is interpreted in topologies according to the above, then projectivity is not preserved. 
We say that $R(U)$ has maximal rank in $\mathcal{D}_{L^{1}}$, if $\dim_{reg} R(U)=\dim R(U)$.

\newtheorem{lemma10}[def1]{Proposition (Maximal rank)}
\begin{lemma10}
 Assume $d U(f)=0$ implies $(I-U) f \in C^{\infty}$. Then $U$ has maximal rank ($ U \neq 0$) over $\mathcal{D}_{L^{1}}$.
\end{lemma10}

Assume $A=\{ f \in \mathcal{D}_{L^{1}} \quad d U(f)=0 \}$ and $reg A=\{ Uf \in C^{\infty} \}$.
Maximal rank means that $\dim A = \dim reg A$. Given U projective, the result is clear. 
When $(I-U)$ corresponds to a reduced measure $d V$ with $d V \bot d U$, we have  maximal rank. 
Consider $E$ very regular, such that $XE=I$ with $\mbox{ ker }E=\{ 0 \}$, where $d U=X dt$, then through the conditions, $UE$ is 
very regular and $E$ is projective on $\mathcal{D}_{L^{1}}$.
Note that point wise topology means $(I- U) f^{2} \sim ((I-U) f)^{2}$, that is the corresponding 
zero set is algebraic.

Assume $U^{\bot}f =0$ (1-polar) implies $f \in C^{\infty}$. In particular $U^{\bot}(d f)=d U^{\bot}(f)=0$
implies $f \in C^{\infty}$. Thus, $d U^{\bot}$ is homogeneously hypoelliptic. Given $U$ projective, 
we have that $d U$ is regularizing and we have maximal rank. Assume $d U + d U^{\bot} - d I= d V$, 
such that $d I - d V \rightarrow d I + d V$ is projective. Given $d V^{\bot}$ hypoelliptic, with 
$d V^{\bot} \bot d V$, then d U has maximal rank. When the dimension for $d V^{\bot}$ is the co dimension 
to isolated singularities on the support for d V, it is sufficient to prove $\mbox{ ker }V^{\bot} \subset L^{1}$.

\newtheorem{lemma101}[def1]{Proposition (Main result)}
\begin{lemma101}
 The spiral has maximal rank over $\Gamma$.
\end{lemma101}

Assume $d V=\alpha d U_{1}$. Consider $(x,y) \rightarrow (u_{1},u_{2})$ and 
$\Gamma=\{ \alpha(u_{1},u_{2})=\alpha(u_{2},u_{1}) \}$. Then we have existence of $\beta \leq 1$, 
with $\Delta \beta=0$ over $\Gamma$ and $\beta=1$ on $\Gamma$. Assume $\Omega$ a neighborhood 
$\Gamma$ under $(u_{1},u_{2}) \rightarrow (x,y)$, so 
that $\beta \geq 0$ on $\Omega$. Then $d V$ has maximal rank on $\Gamma$.

\subsection{The cylinder}
Assume $d U$ BV with support on a cylinder, for instance $d U=\alpha d U_{1}$,
$d U = \beta d U_{2}$, with $\alpha,\beta$ regular. 
When $\Omega$ is a regel domain, then $\log \Omega$ is a translations domain, that is $\Omega$ 
is of order 0 and $(\Omega,\Omega^{\bot})$ defines a cylinder. We assume existence of 
$U_{1} \rightarrow U_{2}$ through classical theory. Further, we assume existence of a neighborhood of 
the boundary $\Gamma=\{ U^{\bot}=U \}$ that is cylindrical.

$L^{2}$ is not nuclear, that is $I : L^{2} \rightarrow L^{2}$, but $I \notin L^{2}$. In $L^{2}$
we can, given $W$ a normal movement, motivate that $W \simeq W^{\bot}$. Assume 
$R(W) \bigoplus A=L^{1}$, where A is finite dimensional. \emph{Assume $d W \in (I)'$ and $d W \rightarrow d I$,
we then do not have that $d I \in (I)'$, if $(I)'$ is not nuclear}.

Assume $d U$ reduced to $d U_{1}$, such that  $U$ intersects all leaves $U^{\bot} =0$ and all 
leaves in the same manner, that is $(d U)^{\bot}$ is defined by $d V$, without $d V$ changing 
character (simply connected). Note that given $\{ d V=0 \}=\{ 0 \}$, $U$ must be projective. 
Given $\Omega$ a translation domain, we have that $\Omega \neq \Omega^{*}$, that is 
\emph{$d I \notin (I)'$. Thus, $d I \in (I)'$ means that $\Omega$ is not a translations domain }, 
that is $d V=\rho d I$, where $\rho$ has constant surfaces. Assume critical points principally 
defined and $\{ d (V - \rho I)=(V - \rho I)=0 \}$ a removable set. Further, $d (V - \rho I)$ 
locally  BV in a neighborhood of this set. 

\newtheorem{prop11}[def1]{Lemma (The approximation property)}
\begin{prop11}
When $(I)' \ni d W \rightarrow d I$ in weak topology and $(I)'$ is nuclear, then $d I \in (I)'$.  
When W is projective and $d W^{\bot}$ a reduced measure, the limit is regular. 
\end{prop11}

$d W \rightarrow d I$ with $(I)'$ nuclear, implies $d I \in (I)'$ through the approximation 
property, that is $d W=\alpha d I$ with $\alpha \rightarrow 1$ regularly. We assume $I_{d W}=<d W,\varphi>=<W, d \varphi>$, with 
uniform convergence. Given $R(W)$ finite dimensional, we have that $W$ reflexive is sufficient 
for nuclearity. When $d W$ is not projective,
we consider $d W + d W^{\bot} - I=d V$. Assume in particular that $(d V)^{\bot}$ is reduced, 
then $d W$ must be projective. Define $\Gamma=\{ d V = d V^{\bot} \}$
and $\Gamma_{0}=\{ \gamma \in \Gamma \quad d V=0 \}$. When $d V \varphi=0$ implies $\parallel \varphi \parallel=0$,
we have that $\widehat{\varphi}=0$, that is $d V$ is reduced and $\Gamma_{0}=\{ 0 \}$. Given $d V$ reduced, 
we have that $\Gamma=\Gamma_{0}$.

\subsection{Nuclearity}
 
 Assume $R(U)^{\bot} \simeq X_{0}$ and $\Phi$ the projection $R(U) \rightarrow X_{0}$. Assume further
 $R(U)$ is defined by $d U$ BV and $X_{0}$ by $d U^{\bot}$ BV. Given $\Phi$ corresponds to 
 regularizing action, U preserves hypoellipticity. Thus, we have that $d U, dU^{\bot}$ do not both 
 preserves hypoellipticity. Assume $\Phi : \mathcal{G}_{HE} \rightarrow \mathcal{G}_{HE}^{\bot}$ 
 preserves hypoellipticity, then $\mathcal{G}_{HE}$ does not have the approximation property.
 
 \newtheorem{lemma3}[def1]{Lemma ($\mathcal{G}$ preserving HE is not nuclear)}
 \begin{lemma3}
  Assume $\Phi : U \rightarrow U^{\bot}$, where $U \in \mathcal{G}$, then the spiral is
  represented by $\Phi = Id$.
 \end{lemma3}
 
 Assume $d U(f) =\{ G,f \}$, then iteration is given by $\{ G, \{ G,f \} \}$.  We claim that 
 the spiral gives a non-closed extension (non-algebraic). Given ${}^{t} U$ corresponds to $U^{\bot}$
 with $U \in \mathcal{D}_{L^{1}}'$, then assuming a very regular boundary, $U \rightarrow {}^{t} U$ corresponds to an algebraic 
 continuation  (modulo $C^{\infty}$). \emph{However, given $U^{2} f=Uf + Vf$, for some $V \neq 0$, 
 then $U+V$ is not necessarily algebraic}. We assume that in a neighborhood of the boundary $U,{}^{t} U$ preserve rank. 
 Assume $(U-I) U f = V f$, we then have $U f= (U-I)^{-1}V f$. When $U^{2},V$ are closed, 
 we do not have that U is closed. Note that $N(V)=N(U) \cup N(U-I)$ (cf. very regular boundary).
 Note $U \in \mathcal{E}'$ does not imply $(U-I)^{-1} \in \mathcal{E}'$, that is $\mathcal{E}'$
 is a discontinuous convolution algebra. 

\bibliographystyle{amsplain}
\bibliography{sos}

\end{document}